\newif\ifstandardtemplate\standardtemplatetrue%
\newif\ifelseviertemplate\elseviertemplatefalse%
\newif\ifspringertemplate\springertemplatefalse%
\newif\ifwileytemplate\wileytemplatefalse%
\newcommand{\mytitle}{Embedding structures in continua:\\ linear models and finite element discretizations}
\newcommand{\myabstract}{This work describes models and numerical approximations that
  describe the mechanical behavior of deformable continua with embedded structural members, such as
  rigid bodies, beams, shells, etc. The continuum formulation extends an idea first presented in the
  context of the Arlequin method and constrains the kinematics of the two types of bodies to be
  compatible in the energy sense. In the article, we exploit the shared similarities of all
  structural theories to introduce a general framework for energetically coupling the latter with
  continua. In addition, we show that the problems, as well as their finite element approximations,
  are well-posed. Numerical examples of bodies with inclusions, fibers, and embedded surfaces are
  provided to illustrate the generality and robustness of the approach.}

\newcommand{\mypackages}{%
  \usepackage[utf8]{inputenc}
  \usepackage{amssymb}
  \usepackage{amsmath}
  \usepackage{amsthm}
  \usepackage{graphicx}
  \usepackage{caption}
  \usepackage{subcaption}
  \usepackage{natbib}
  \usepackage{xfrac}
  \usepackage{todonotes}
  \usepackage{siunitx}
  \usepackage{changes}
  \usepackage{comment}
  \definechangesauthor[name=Ignacio, color=blue]{IRO}
  \definechangesauthor[name=Ignacio, color=orange]{DPG}
  \usepackage{cleveref}
  \graphicspath{{./}{Figures/}{figures/}}
  \usepackage{enumitem}
  \setlist[enumerate,1]{label = \emph{\alph*}),
                        ref   = \theenumi.\emph{\alph*}}
}

\newcommand{\mymacros}{%
  \newcommand{\concept}[1]{\textbf{\emph{##1}}}
  \newcommand{\combo}[2]{{\{##1,##2\}}}
  \newcommand{\defined}{:=}
  \newcommand{\dev}{{\mathop{\mathrm{dev}}}}
  \renewcommand{\div}{{\mathop{\mathrm{div}}}}
  \newcommand{\mbs}[1]{{##1}}
  \newcommand{\pairing}[2]{\langle{##1},{##2}\rangle}
  \newcommand{\dd}[2]{\frac{\mathrm{d} ##1}{\mathrm{d} ##2}}
  \newcommand{\pd}[2]{\frac{\partial{##1}}{\partial{##2}}}
  \newcommand{\set}[1]{\left\{##1\right\}}
  \newcommand{\structinner}[2]{\langle\kern-1pt\langle ##1, ##2 \rangle\kern-1pt\rangle}
  \newcommand{\trace}{{\mathop{\mathrm{tr}}}}
  \newcommand{\triplenorm}[1]{\|\kern-1pt|##1|\kern-1pt\|}
  \newcommand{\uptohere}{\centerline{\textcolor{blue}{\rule{6cm}{0.2cm}}}}
  \let\oldLambda=\Lambda\renewcommand{\Lambda}{\mathit{\oldLambda}}
  \let\oldGamma=\Gamma\renewcommand{\Gamma}{\mathit{\oldGamma}}
}

\newcommand{\base}{\mathcal{C}}
\newcommand{\fiber}{\mathcal{F}}

\newcommand{\dV}{\mathrm{d}V}

\newcommand{\hs}{h_{\Omega}}
\newcommand{\hb}{h_{\base}}
\newcommand{\uhb}{\bar{u}_{\hb}}
\newcommand{\vhb}{\bar{v}_{\hb}}
\ifstandardtemplate%
\documentclass[10pt,a4paper]{article}
\mypackages%
\usepackage[noblocks]{authblk}
\mymacros%

\newcommand{\mybibstyle}{unsrt}

\newtheorem{theorem}    {Theorem}[section]

\theoremstyle{definition}

\newtheorem{examplex}[theorem]{$\triangleright\;$Ejemplo}

\theoremstyle{remark}

\title{\mytitle}
\author[1]{David Portillo}
\author[1,2]{Ignacio Romero}

\affil[1]{Dept. of Mechanical Engineering, Universidad Politécnica de Madrid,
  Jos\'{e} Guti\'{e}rrez Abascal, 2, 28006 Madrid, Spain}
\affil[2]{IMDEA Materials Institute, Eric Kandel 2, 28096 Getafe, Madrid, Spain}

\begin{document}
\maketitle
\begin{abstract}
  \myabstract
\end{abstract}
\fi

\ifspringertemplate%
\documentclass[smallcondensed]{svjour3}
\mypackages%
\mymacros%

\title{\mytitle%
\thanks{the thanks}}
\author{\ldots \and Ignacio Romero \and \ldots}
\journalname{The journal name}

\institute{I. Romero \at%
  IMDEA Materials Institute, Eric Kandel 2, Tecnogetafe, Madrid 28906, Spain\\
  Universidad Polit\'ecnica de Madrid, Jos\'e Guti\'errez Abascal, 2, Madrid 29006, Spain\\
  \email{ignacio.romero@imdea.org}
  \and
  XX \at%
  XXX\\
  \email{XXX}}

\titlerunning{\mytitile}
\authorrunning{I. Romero}

\date{Received: date / Accepted: date}

\begin{document}
\maketitle
\begin{abstract}
  \myabstract%
\end{abstract}
\fi

\ifelseviertemplate%
\documentclass[preprint,11pt]{elsarticle}
\mypackages
\mymacros
\newcommand{\mybibstyle}{unsrt}
\newtheorem{theorem}{Theorem}

\journal{Computer Methods in Applied Mechanics and Engineering}
\begin{document}
\begin{frontmatter}

\title{\mytitle}

\author[1]{David Portillo}
\ead{david.portillo@upm.es}

\author[1,2]{Ignacio Romero\corref{cor1}}
\ead{ignacio.romero@upm.es}

\address[1]{Universidad Polit\'ecnica de Madrid, Spain}
\address[2]{IMDEA Materials Institute, Spain}

\cortext[cor1]{Corresponding author. ETS Ingenieros Industriales,
  Jos\'e Guti\'errez Abascal, 2, Madrid 28006, Spain}

\begin{abstract}
\myabstract%
\end{abstract}

\end{frontmatter}

\fi

\ifwileytemplate%
\documentclass[doublespace,times]{nmeauth}
\mypackages%
\mymacros%
\usepackage{natbib}
\message{Compiling with Wiley template}

\begin{document}
\runningheads{I. Romero}{\dots}

\title{\mytitle}
\author{Ignacio Romero\affil{1}\affil{2}\corrauth}

\address{\affilnum{1} ETSII, Universidad Politécnica de Madrid,
         Jos\'{e} Guti\'{e}rrez Abascal, 2, 28006 Madrid, Spain\break%
         \affilnum{2} IMDEA Materials Institute, Eric Kandel 2, 28096 Getafe, Madrid, Spain}

\corraddr{Dpto.~de Ingenier\'{\i}a Mec\'{a}nica; E.T.S. Ingenieros Industriales;
José Gutiérrez Abascal, 2; 28006 Madrid; Spain. Fax (+34) 91 336 3004}

\begin{abstract}
\end{abstract}

\maketitle
\fi




\section{Introduction}
\label{sec-intro}
In solid mechanics, structural models such as beams, plates, and shells are employed to analyze the
behavior of bodies with special geometries and simple kinematics. When using numerical
approximations, the latter employ fewer degrees of freedom that their continuum counterparts and,
when possible, they are favored by analysts. Often, however, the structures under study do not have
any special geometry and thus their behavior can only be studied by solving a problem of classical
solid mechanics of the type commonly considered in monographs and textbooks on continuum mechanics.

Models of structural members and continuum bodies are both well-known. Their numerical
discretization, for example in the context of finite elements, is also well-known and nowadays
robust and accurate formulations are widespread. It remains, however, one aspect of these problems
that remains open: the description of mechanical problems that combine structural and continuum
bodies that must co-exist harmonically, where loads and boundary conditions could be applied on
either type and where internal forces should be transferred among bodies while preserving the
conservation laws that govern equilibrium.

This type of mechanical problems with governing equations for bodies with mixed dimensionality is
helpful to describe numerous situations in Engineering and Materials Science. For instance, a
concrete slab with steel reinforcing bars can be conveniently described using a continuum domain for
the slab, structural beams for the rebars, and their interaction. In the context of biomaterials,
tendons provide a clear example of reinforcing of a matrix with strong (collagen) fibers, similarly
to what occurs in plants such as flax or bamboo.

In this work we are interested in analyzing the mechanical behavior of bodies consisting of a matrix
with embedded structures of dimensions 0, 1, or 2. From the mathematical point of view, the difficulty
of such task is due to the fact that each of the involved problems is governed by equations that are
different in nature. Strictly speaking, the functional spaces where the solution of the continuum
bodies live are Hilbert spaces that, in general, do not have well-defined traces in
lower-dimensional manifolds, like curves or surfaces \cite{angelo2008ca}. Hence, ``connecting'' the solution fields in
the continua and the structures requires special glueing techniques that will have to be developed
\emph{ad hoc}.

The relevance of mixed-dimensionality problems has motivated the appearance of different avenues
aimed ad modeling and/or finding numerical solutions to this class of coupled problems. When focusing on the formulation of numerical
methods, interface finite elements have been used for several years
\cite{bournival2010id,klarmann2022xz}. Also, constrained problems have been proposed that
enforce the kinematic compatibility or equilibrium of force resultants
\cite{shim2002zi,mccune2000kk,romero2018iu,romero2023jg,schenk2023rs}. Moreover, using Nitsche's method, it is
possible to enforce the weak compatibility of structures and continua
\cite{nguyen2013ws,hansbo2022fz}, but also using mortar finite
elements (e.g. \cite{steinbrecher2020cm,Steinbrecher2022a,firmbach2023up}). Finally, some authors
have noted that it is more natural to couple structures and three-dimensional bodies when the latter are modeled
using Cosserat continua \cite{sky2024qd}.

In this work, we propose a glueing procedure that stems from the Arlequin method \cite{dhia2001iu,dhia2005iu},
a technique developed to simplify the connection of (continuum) bodies. Vaguely speaking, this
method imposes the compatibility of the displacement fields of the interacting bodies in the energy
sense, ensuring the well-posedness of the joint formulation. In the current work, we extend this
idea to formulate linking formulations that guarantee compatible displacement fields of the
deformable matrix and the embedded structures. As in the Arlequin method, via a suitable definition
of the interface terms, the joint equilibrium problem will be shown to be well-posed. Remarkably,
straightforward mixed finite elements can be shown to be stable.

The techniques introduced in the current work serve to glue rigid inclusions, line and surface
reinforcements to two- and three-dimensional bodies. In the article, we will put forward a general
framework that describe all these problems in a unified fashion. Additionally, it will be possible
to prove in a compact fashion the well-posedness of a large class of embedding problems.

The numerical examples that we consider showcase different types of embeddings, namely, beam and
shells in solids. However, the generality of the framework indicates that other linking problems can
be adapted by using simple modifications.

Let us mention, in closing, that the boundary value problems considered in this work are all
mechanical. However, the translation of the ideas presented to other equilibrium problems (for
example, thermal) are fairly straightforward.

The remainder of the article has the following structure. In Section~\ref{sec-model}, we provide a succinct
summary of the governing equations of continuum bodies and structures, all in the framework of small
strain elasticity. This is a vast subject and we only intend to highlight the common structure of
all structural bodies and introduce the required notation. In Section~\ref{sec-embedding} we
describe the coupling problem. The continuum equations of linked bodies are introduced and their
well-posedness, discussed. Then, their finite element discretization is addressed in
Section~\ref{sec-fem}. While the discretization is completely standard, it remains to show that the
finite dimensional projection leads to discrete problems that are well-posed. In general, this task
is not easy but it will be shown that the framework previously introduced simplifies the proofs.
Section~\ref{sec-examples} illustrates the proposed formulations using examples with varied
geometries. Finally, a summary of the main results in the article and an outlook towards future
extensions are provided in Section~\ref{sec-summary}. Some of the more technical results are left
for Appendix \ref{app-instabilities}.


\section{Linear models for continuum bodies and structures}%
\label{sec-model}
This article aims at formulating models and numerical methods that connect, in a precise
manner to be defined later, bodies described with mathematical models that are
intrinsically different. More precisely, we are interested in coupling \emph{continuum}
with \emph{structural} models, all in the framework of mechanics.

There is a large collection of models in solid and structural mechanics, and the methods introduced in this work are general enough to connect many of them. Given the generality of the goal and the multiplicity of possible combinations, we start by reviewing succinctly several well-known \emph{linear} models. By presenting them in a unified way, it will be simpler to introduce later the connecting equations.

\subsection{Linearized elasticity}%
\label{subs-linear-elas}
We first consider linearized mechanics restricting, for simplicity, our exposition to
three-dimensional elasticity. This is the template for more complex continuum models and all exploit the relative simplicity of the topology and geometry of Euclidean space. Given the isotropy of Euclidean space, defining (weak) derivatives of the relevant fields is straightforward. Also, the natural energy norms and inner products are directly derived from standard counterparts for Hilbert spaces.

The model of linearized elasticity is the simplest one that can be used to represent the behavior of continuum bodies and, as a result, it is well-understood. Numerous monographs describe the problem in detail from all points of view. For the sequel, it suffices to define elastic bodies as bounded, open, subsets $\Omega\subset \mathbb{R}^3$ with smooth boundary and points denoted as $\mbs{x}$.

The deformation of these bodies is customarily described via the displacement field $\mbs{u}:\Omega\to \mathbb{R}^3$. This vector field must belong to $X \defined [H^1_{\partial}(\Omega)]^3$, the Hilbert space of square-integrable functions with square-integrable derivatives, and vanish in $\partial_D\Omega$, a subset of the boundary of $\Omega$ that has a nonzero measure. If $\ell$ is a characteristic length of the domain, for example, its diameter, two vector fields $\mbs{u},\mbs{v}\in X$ have the inner product
\begin{equation}
  \label{eq-inner-linearized}
  \langle \mbs{u} , \mbs{v} \rangle_{\Omega}
  \defined
  \int_{\Omega}
  \left(
    \mbs{u}(\mbs{x})\cdot \mbs{v}(\mbs{x})
    +
    \ell^2\;
    D\mbs{u}(\mbs{x})\cdot D\mbs{v}(\mbs{x})
  \right)
  \; \mathrm{d} V\ ,
\end{equation}
where the dot indicates the contraction of all indices, we have employed the notation
\begin{equation}
  \label{eq-D}
  D\mbs{u}(\mbs{x})
  \defined
  \pd{u_i}{x_j}\; \mbs{E}_i\otimes \mbs{E}_j\ ,
\end{equation}
and $\{\mbs{E}_i\}_{i=1}^3$ is the Cartesian basis of $\mathbb{R}^3$. Let us note that, due to
Korn's inequality \cite{ciarlet1988ux}, the norm induced by this inner product is equivalent to the
\emph{energy} norm defined over displacement fields.

A key feature of this mechanical model --- one that will be exploited repeatedly -- is that it possesses a variational statement. Given some external volume forces and surface tractions, there
exists a potential energy functional $V_{LE}:X\to \mathbb{R}$ such that the
equilibrium displacement of the solid is the one that solves
\begin{equation}
  \label{eq-le-equilibrium}
  \inf_{\mbs{u}\in E} V_{LE}(\mbs{u})\ .
\end{equation}
In linear problems the potential $V_{LE}$ is quadratic and there exists a linear form $a_{LE}:X\times X\to \mathbb{R}$ and a linear form $f_{LE}:X\to \mathbb{R}$ such that
\begin{equation}
  \label{eq-le-weak}
  V_{LE}(\mbs{u}) \defined
  \frac{1}{2} a_{LE}(\mbs{u},\mbs{u}) - f_{LE}(\mbs{u})\ .
\end{equation}
The properties of these two forms will play a key role in the well-posedness of the models and numerical methods later introduced.

\subsection{Abstract linear structural models}%
\label{subs-reduced}

As advanced, there are several reduced dimensional theories (beams, shells, membranes, \ldots).
All can be understood as truncated asymptotic expansions of the equations that govern the continuum
theory. Here, we focus on a class of linear models and give them a common description that will be
exploited in Section~\ref{sec-embedding}. Let us note that not all structural models can be framed
in the formalism that follows.

A mechanical structure is a deformable body that can be identified, just as in
Section~\ref{subs-linear-elas}, with a bounded domain $\mathcal{S}\subset \mathbb{R}^3$. However, and in
contrast with general bodies, a structure can also be identified with a \emph{fiber bundle} \cite{epstein2010vi}. Hence, there
exists two manifolds, the base $\mathcal{C}$ and the fiber $\mathcal{F}$, and two surjective projections $\Pi_{\mathcal{C}}:\mathcal{S}\to
\mathcal{C}$
and $\Pi_{\mathcal{F}}:\mathcal{S}\to \mathcal{F}$. The manifolds $\mathcal{C},\mathcal{F}$ have
dimension $n\ge0$ and $m>0$, respectively, with $n+m=3$; moreover, the fiber can be identified with a bounded subset
$\mathbb{R}^{m}$. All structural models possess a characteristic dimension $\ell$. It can be the radius of gyration
of a rigid body, the thickness of a plate, etc.

For convenience, we define $\Pi \defined (\Pi_{\mathcal{C}},\Pi_{\mathcal{F}})$ and for all particles $\mbs{x}\in \mathcal{S}$ we write
\begin{equation}
  \label{eq-split}
  \Pi(\mbs{x})
  =
  (\mbs{\sigma},\mbs{\xi}) = (\Pi_{\mathcal{C}}(\mbs{x}), \Pi_{\mathcal{F}}(\mbs{x}))
  \ .
\end{equation}
We call $\mbs{\sigma}$ and $\mbs{\xi}$ the base and fiber coordinates, respectively, of a point
$\mbs{x}\in \mathcal{S}$. The inverse
of the map $\Pi$ is the embedding $\mathcal{E}:\mathcal{C}\times \mathcal{F}\to \mathcal{S}\subset
\mathbb{R}^3$. Thus, for any integrable function $f:\mathcal{S}\to \mathbb{R}$, we have that
\begin{equation}
  \label{eq-jacobian1}
  \int_{\mathcal{S}} f(\mbs{x}) \; \mathrm{d} V
  =
  \int_{\mathcal{C}}\int_{\mathcal{F}} (f\circ \mathcal{E})
  (\mbs{\sigma},\mbs{\xi})
  \; J(\mbs{\sigma},\mbs{\xi})
  \; \mathrm{d} \mathcal{F} \; \mathrm{d} \mathcal{C}\ ,
\end{equation}
where $J$ is the Jacobian of the embedding defined as
\begin{equation}
  \label{eq-jacobian}
  J(\mbs{\sigma},\mbs{\xi})
  \defined
  \left|
    \pd{\mathcal{E}(\mbs{\sigma},\mbs{\xi})}{(\mbs{\sigma},\mbs{\xi})}\right|\ .
\end{equation}

To introduce the kinematics of a mechanical structure we require two smooth mappings
$\mbs{\Sigma}:\mathcal{C}\to \mathbb{R}^{3}$ and $\mbs{\Theta}:\mathcal{C}\to
\mathcal{L}(\mathbb{R}^{m},\mathbb{R}^3)$. Here, $\mathcal{L}(\mathbb{R}^{m},\mathbb{R}^3)$ refers to the vector space of linear maps from the fiber to $\mathbb{R}^{3}$. The displacement of the structure is the mapping $\Psi\combo{\Sigma}{\Theta}:\mathcal{C}\times \mathcal{F}\to \mathbb{R}^3$ of the form
\begin{equation}
  \label{eq-struct-psi}
  \Psi\combo{\Sigma}{\Theta}(\mbs{\sigma},\mbs{\xi})
  \defined
  \mbs{\Sigma}(\mbs{\sigma}) +
  \mbs{\Theta}(\mbs{\sigma})*\mbs{\xi}\ ,
\end{equation}
where $(\mbs{\sigma},\mbs{\xi})$ are, respectively, the coordinates of $\mathcal{C}$ and
$\mathcal{F}$ and $\mbs{\Sigma},\mbs{\Theta}$ define the displacement of the base and fibers, respectively. Let us note that the map~$\mbs{\Theta}$ depends on the base coordinates
and acts \emph{linearly} on the fiber coordinates. This operation, not precisely defined yet, is indicated with the symbol `$*$' and depends on the structural model.

A field of the form~\eqref{eq-struct-psi} is univocally described by the two generalized
displacements $\Sigma$ and $\Theta$. Then, slightly abusing the notation, we refer to the function
$\Psi\combo{\Sigma}{\Theta}$ simply with the bracketed pair $\combo{\Sigma}{\Theta}$. These displacements
must belong to the space
\begin{equation}
  \label{eq-struct-space}
  Y
  \defined
  \left\{ \combo{\Sigma}{\Theta},\
    \mbs{\Sigma} \in [H^1(\mathcal{C})]^3, \
    \mbs{\Theta}\in [H^1(\mathcal{C})]^{r}
  \right\}
  \ ,
\end{equation}
where $0\le r\le 3$ is the dimension of the range of $\Theta$, and depends on the structural model,
more precisely to the allowed motions of the fiber. The extreme case $r=0$ corresponds to structural
models that do not contemplate fiber rotations (e.g., bars and membranes).


The bundle structure of a mechanical model allows us to express its displacements and other fields
 either as maps from $\mathcal{S}\subset\Omega$ or from $\mathcal{C}\times
\mathcal{F}$. In particular, if $\mbs{u},\mbs{v}$ are two vector fields defined on $\mathcal{S}$, we can
write their inner product $\langle \mbs{u}, \mbs{v}\rangle_{\mathcal{S}}$ using
expression~\eqref{eq-inner-linearized}. Crucially for our developments, we can define an alternative
inner product for displacement fields in the structure. To see this, let $U \defined u\circ \mathcal{E}$ and
$V \defined v\circ \mathcal{E}$ and define
\begin{equation}
  \label{eq-struct-inner2}
  \structinner{\mbs{U}}{\mbs{V}}
    \defined
      \int_{\mathcal{C}}\int_{\mathcal{F}}
      \left(
      \mbs{U}(\mbs{\sigma},\mbs{\xi})\cdot \mbs{V}(\mbs{\sigma},\mbs{\xi})
      +
      \ell^2\;
      D\mbs{U}(\mbs{\sigma},\mbs{\xi})\cdot D\mbs{V}(\mbs{\sigma},\mbs{\xi})
      \right)
\; \mathrm{d} \mathcal{F} \; \mathrm{d} \mathcal{C}
      \ .
\end{equation}
Note that, in general, $\langle \mbs{U},\mbs{V}\rangle_{\mathcal{S}} \neq
\structinner{\mbs{U}}{\mbs{V}}$ because the inner product~\eqref{eq-struct-inner2}
does not contain the Jacobian~\eqref{eq-jacobian} of the embedding. This is a common approximation in structural theories, one that simplifies their equations and, often, allows
for analytical solutions.

As in the case of the continuum solid, a key feature of all structural theories is that they also admit a variational principle, one that can be obtained from Eq.~\eqref{eq-le-equilibrium} using the appropriate simplifications. Thus, there exists a potential energy $V_{ST}$ that depends on the deformation of the structure as well as the external loads such that
\begin{equation}
  \label{eq-inf-structure}
  \combo{\Sigma}{\Theta}
  = \arg\inf_{\combo{\Upsilon}{\beta}\in Y} V_{ST}(\Upsilon,\beta)\ .
\end{equation}
The form of the potential energy $V_{ST}$ depends, of course, on the specific structural model. However, in all cases, the linearity of the model requires that the potential energy be quadratic and thus there must exist a bilinear form $a_{ST}:Y\times Y\to \mathbb{R}$ and a linear form $f_{ST}:Y\to \mathbb{R}$ such that
\begin{equation}
  \label{eq-struct-forms}
  V_{ST}(\Upsilon,\beta)
  \defined
  \frac{1}{2} a_{ST}(\Upsilon,\beta;\Upsilon,\beta) - f_{ST}(\Upsilon,\beta)\ .
\end{equation}

A wide range of structural models, including shells, beams, bars, and rigid bodies, can be represented in an abstract manner using the previously introduced framework. This representation can be expressed using the following bracketed pair $ \combo{\Sigma}{\Theta}(\mbs{\sigma},\xi)$,
\begin{equation}
  \label{eq-general-psi}
  \combo{\Sigma}{\Theta}(\mbs{\sigma},\xi) =
  \mbs{\Sigma}(\mbs{\sigma}) + \mbs{\Theta}(\mbs{\sigma}) * \Xi(\xi) \ ,
\end{equation}
where $\Xi = \xi_i E_i \in \mathcal{F}$ is a vector in the fiber and $\mbs{\Theta}(\mbs{\sigma})$ is a tensor that
represents an infinitesimal rotation of the fiber, thus $\mbs{\Theta}:\mathcal{C}\to so(3)$.

The inner product~\eqref{eq-struct-inner2} can then be developed from the derivative of the
displacement, yielding
\begin{equation}
  \label{eq-general-grad}
  D\combo{\Sigma}{\Theta}(\mbs{\sigma},\xi) =
  \pd{\combo{\Sigma}{\Theta}(\mbs{\sigma},\xi)}{\sigma_\alpha}\otimes \mbs{E}^{\alpha} +
  \pd{\combo{\Sigma}{\Theta}(\mbs{\sigma},\xi)}{\xi^i}\otimes \mbs{E}^i\ ,
\end{equation}
where $\alpha = 1,..,n$ and $i=1,...,m$. Hence, the structural inner
product~\eqref{eq-struct-inner2} becomes
\begin{equation}
  \label{eq-general-inner}
  \begin{split}
      \structinner{\combo{\Sigma}{\Theta}}{\combo{\Upsilon}{\Lambda}}
    =& \int_{\mathcal{C}}
       |\mathcal{F}|
    \left(
       \mbs{\Sigma}(\sigma)\cdot \mbs{\Upsilon}(\sigma) + \ell^2\, \mbs{\Sigma}'(\sigma)\cdot \mbs{\Upsilon}'(\sigma)
       \right) \; \mathrm{d}\mathcal{C}
       \\
    &+
\int_{\mathcal{C}}
    \left(
      \mbs{\theta}(\sigma)\cdot \mbs{i}_{\mathcal{F}}\,\mbs{\lambda}(\sigma) +
      \ell^2\, \mbs{\theta}'(\sigma)\cdot \mbs{i}_{\mathcal{F}}\, \mbs{\lambda}'(\sigma)
      \right) \;\mathrm{d}\mathcal{C}
    \\
    &+
    \int_{\mathcal{C}}
    |\mathcal{F}|
      \ell^2\;
      \mbs{\theta}(\sigma)\cdot
      (\mbs{I} + \mbs{E}^i\otimes \mbs{E}^{i})\mbs{\lambda}(\sigma)
      \; \mathrm{d}\mathcal{C}
    \ .
  \end{split}
\end{equation}
Here, $\lambda$ is the vector field of axial vectors of the skew-symmetric tensor field $\Lambda$.
Note that if $\dim(\mathcal{C})=0$, as for rigid solids, the integral should be understood as the
integrand evaluated at the only point of the base manifold.

The common formulation of structural models simplifies the analysis and presentation of tied formulations. For concreteness, we work out some examples next. For them, we specify the base and fiber manifolds and the restricted expressions for their displacement fields. The examples selected have, respectively, base manifolds of dimensions zero, one, and two, illustrating the generality of the approach.

\subsubsection{A linearized rigid body}
The simplest generalized structural model is afforded by a rigid body in the linearized theory. This
is a region  $\mathcal{S}\subset\mathbb{R}^3$ that can only undergo infinitesimal rigid body motions, that is,
translations and infinitesimal rotations.

The base manifold for this body is a single point and, without loss of generality, we can select it
to be its center of mass whose position is indicated as~$\mbs{\sigma}$; the fiber is the body itself, a bounded set in $\mathbb{R}^3$. The dimensions of the base and fiber are, therefore, zero and three, respectively. Given a point $\mbs{x}\in \mathcal{S}$ we have
\begin{equation}
  \label{eq-lrb-projections}
  \Pi_{\mathcal{C}}(\mbs{x}) = \mbs{\sigma}\ ,
  \qquad
  \Pi_{\mathcal{F}}(\mbs{x}) = \mbs{x}-\mbs{\sigma}\ .
\end{equation}
The deformations of this kind of solids are described by the  map
$\mbs{\Sigma}:\mathcal{C}\to \mathbb{R}^3$ and the linear map $\mbs{\Theta}:\mathcal{C}\to so(3)$,
the set of skew-symmetric tensors. Since the base manifold consists of a single point, we write
$\Sigma,\Theta$ instead of $\Sigma(\sigma)$ and $\Theta(\sigma)$.
Combining these two mappings as indicated in Eq.~\eqref{eq-struct-psi}, we find the displacement
\begin{equation}
  \label{eq-lrb-motion}
  \combo{\Sigma}{\Theta}(\sigma,\mbs{\xi})
  \defined
  \Sigma + {\mbs{\Theta}}\, \mbs{\xi}\ .
\end{equation}
Let us recall that the multiplication of a skew-symmetric tensor ${\mbs{\Theta}}$ over a vector
$\mbs{\xi}$ is equivalent to the vector product $\mbs{\theta}\times \mbs{\xi}$, where $\mbs{\theta}$
is the axial vector of the tensor $\mbs{\Theta}$. In Eq.~\eqref{eq-lrb-motion} we recognize the composition of a translation and an infinitesimal rotation.

The inner product for linearized rigid solids can be calculated by combining the general expression~\eqref{eq-struct-inner2} with the specific form of the displacement given by Eq.~\eqref{eq-lrb-motion}. For that, consider two rigid displacement fields
\begin{equation}
  \label{eq-solid-inner-example}
  \combo{\Sigma}{\Theta}(\sigma,\mbs{\xi}) = \Sigma + \mbs{\theta}\times\mbs{\xi}\ ,
  \qquad
  \combo{\Upsilon}{\Lambda}(\sigma,\mbs{\xi}) = \Upsilon + \mbs{\lambda}\times\mbs{\xi}\ ,
\end{equation}
where $\theta,\lambda$ are the axial vectors of $\Sigma$ and $\Lambda$, respectively.
Then, the inner product of these two displacement fields is
\begin{equation}
  \label{eq-inner-solid}
  \begin{split}
  \structinner{\combo{\Sigma}{\Theta}}{\combo{\Upsilon}{\Lambda}}
  &=
  \int_{\mathcal{C}} \int_{\mathcal{F}}
  \left[
    (\mbs{\sigma} + \mbs{\theta}\times\mbs{\xi}) \cdot
    (\mbs{\Upsilon} + \mbs{\lambda}\times\mbs{\xi})
    +
    \ell^2\,
    {\mbs{\Theta}}\cdot {\mbs{\Lambda}}
  \right]
    \; \mathrm{d} \mathcal{F}
    \; \mathrm{d} \mathcal{C}
    \\
  &= |\mathcal{S}|\; (\mbs{\sigma}\cdot \mbs{\Upsilon})
    + \mbs{\theta}\cdot
    \left( \mbs{i}_{\mathcal{S}} + 2\,\ell^2\,\mbs{I} \right) \mbs{\lambda}\ ,
    \end{split}
\end{equation}
where $|\mathcal{S}|$ is the volume of the rigid body, and $\mbs{i}_{\mathcal{S}}, \mbs{I}$ are the inertia and identity tensors, respectively.

The potential energy of a rigid body is denoted as $V_{RB}(\mbs{\Sigma},\mbs{\theta})$ and consists only of the contributions of the external forcing since it can not deform.

\subsubsection{The Timoshenko beam}
The Timoshenko beam is a structural model with a one-dimensional base manifold $\mathcal{C}$ that
corresponds with the curve of centroids of its cross sections. This curve can be described with a
map $\mbs{r}:[0, L]\to \mathbb{R}^3$, where~$L$ is the length of the undeformed beam and the
manifold unique coordinate~$\sigma$ can be identified with the arc-length of $\mbs{r}$. The fiber
$\mathcal{F}$ is the cross-section $A\subset \mathbb{R}^2$ assumed, without loss of generality, to
be constant. We denote as $(\sigma,\mbs{\xi})\in [0,L] \times A$ the base and fiber coordinates,
respectively. Assuming, for simplicity, that the beam is straight, we can define an orthonormal
basis $\{\mbs{E}_i\}_{i=1}^3$, with $\mbs{E}_3=\mbs{r}'(\sigma)$ and $\{\mbs{E}_{\alpha}\}_{\alpha=1}^{2}$ being the principal directions of the cross-section. The projection onto the base is
\begin{equation}
  \label{eq-timo-projections1}
  \Pi_{\mathcal{C}}(\mbs{x}) = \arg\min_\sigma | \mbs{r}(\sigma)-\mbs{x} |\ ,
\end{equation}
and the projection onto the fiber is
\begin{equation}
  \label{eq-timo-projections2}
  \qquad
  \Pi_{\mathcal{F}}(\mbs{x}) =
  \left( \mbs{\xi}\cdot \mbs{E}_1 , \mbs{\xi}\cdot \mbs{E}_2 \right)
  \qquad
  \text{with}
  \qquad
  \mbs{\xi} = \mbs{x}-\mbs{r}(\Pi_{\mathcal{C}}(\mbs{x})) \ .
\end{equation}

The displacement field of a Timoshenko beam is described by two mappings. First, the function $\mbs{\Sigma}:\mathcal{C}\to \mathbb{R}^3$ is used to determine the displacement of the centroids in the cross sections; second, the function $\mbs{\Theta}:\mathcal{C}\to so(3)$ defines the infinitesimal rotation of the cross-section. If $\mbs{\theta}$ is the axial vector of $\mbs{\Theta}$, we can write the displacement function of the Timoshenko beam as
\begin{equation}
  \label{eq-timo-psi}
  \combo{\Sigma}{\Theta}(\sigma,\mbs{\xi})
  \defined
  \mbs{\Sigma}(\sigma) + \mbs{\theta}({\sigma})\times \mbs{\xi}\ .
\end{equation}

The inner product of displacements of the Timoshenko beam possesses a simple expression. To obtain it, consider first two displacement fields of the form
\begin{equation}
  \label{eq-timo-inner-example}
  \combo{\Sigma}{\Theta}(\sigma,\mbs{\xi}) = \mbs{\Sigma}(\sigma) + \mbs{\theta}(\sigma)\times\mbs{\xi}\ ,
  \qquad
  \combo{\Upsilon}{\Lambda}(\sigma,\mbs{\xi}) = \mbs{\Upsilon}(\sigma) + \mbs{\lambda}(\sigma)\times\mbs{\xi}\ ,
\end{equation}
where $\theta(\sigma),\lambda(\sigma)$ are the axial vectors of $\Theta(\sigma),\Lambda(\sigma)$, respectively.
To compute their inner product, let us note that the derivative of a displacement map is
\begin{equation}
  \label{eq-timo-D}
  D\combo{\Sigma}{\Theta}(\sigma,\mbs{\xi}) =
  \pd{\combo{\Sigma}{\Theta}(\sigma,\mbs{\xi})}{\xi_{\alpha}}\otimes \mbs{E}_{\alpha} +
  \pd{\combo{\Sigma}{\Theta}(\sigma,\mbs{\xi})}{\sigma}\otimes \mbs{E}_{3}\ ,
\end{equation}
where Greek repeated indices indicate a sum from 1 to~2. With this definition, we can evaluate the
inner product of two Timoshenko displacement fields as
\begin{equation}
  \label{eq-inner-timo}
  \begin{split}
    \structinner{\combo{\Sigma}{\Theta}}{\combo{\Upsilon}{\Lambda}}
    =& \int_{0}^{L}
       |A|
    \left(
       \mbs{\Sigma}(\sigma)\cdot \mbs{\Upsilon}(\sigma) + \ell^2\, \mbs{\Sigma}'(\sigma)\cdot \mbs{\Upsilon}'(\sigma)
       \right) \;\mathrm{d}\sigma
       \\
    &+
\int_{0}^{L}
    \left(
      \mbs{\theta}(\sigma)\cdot \mbs{i}_{\mathcal{F}}\,\mbs{\lambda}(\sigma) +
      \ell^2\, \mbs{\theta}'(\sigma)\cdot \mbs{i}_{\mathcal{F}}\, \mbs{\lambda}'(\sigma)
      \right) \;\mathrm{d}\sigma
    \\
    &+
    \int_{0}^{L}
      |A|\,\ell^2\;
      \mbs{\theta}(\sigma)\cdot
      (\mbs{I} + \mbs{E}_{\alpha}\otimes \mbs{E}_{\alpha})\mbs{\lambda}(\sigma)
      \;\mathrm{d}\sigma
    \ ,
  \end{split}
\end{equation}
where $|A|$ is the cross-section area, $\mbs{i}_{\mathcal{F}}$ is, as before, the inertia of the fiber, and $\mbs{I}$ is the identity tensor. The form of this product seems more involved than the one in Eq.~\eqref{eq-inner-linearized} but, in contrast with the general definition, it only involves a line integral, a clear consequence of the dimensional reduction of the theory.

The potential energy of an elastic Timoshenko beam is a functional $V_T(\mbs{\Sigma},\mbs{\theta})$ that includes the internal energy due to all its deformation modes as well as the contributions from the external loading. The equilibrium displacements $\combo{\Sigma}{\Theta}$ are the minimizers
\begin{equation}
  \combo{\Sigma}{\Theta} =
  \arg\inf_{\combo{\Upsilon}{\Lambda}\in Y} V_T(\Upsilon,\Lambda)\ .
\end{equation}

\subsubsection{The Mindlin shell}
The Mindlin shell is one of the many models available for thin solids. In this case, the base
manifold $\mathcal{C}$ is the mid-surface of the shell itself, a bounded flat surface~$A$ whose
points in $\mathbb{R}^3$ are given by a two-parameter function $\mbs{r}:\mathcal{C}\to
\mathbb{R}^3$. In this model, the coordinates of $\mathcal{C}$ are $\sigma=(\sigma_{1},\sigma_2)\in
\mathbb{R}^2$. The fiber is the segment $[-t/2,t/2]\in \mathbb{R}$, where $t$ is the thickness of the shell. We define a Cartesian basis $\{\mbs{E}_i\}_{i=1}^3$ with $\mbs{E}_3$ being the unit normal to the plane. Then, the projection maps for a point $\mbs{x}\in \mathcal{S}$ are
\begin{equation}
  \label{eq-mindlins-reisness}
  \Pi_{\mathcal{C}}(\mbs{x}) = \arg\min_{\sigma_1,\sigma_2} |\mbs{x}-\mbs{r}(\sigma_1,\sigma_2)|\ ,
  \qquad
  \Pi_{\mathcal{F}}(\mbs{x}) = (\mbs{x}-\mbs{r}(\Pi_{\mathcal{C}}(\mbs{x})))\cdot \mbs{E}_3\ .
\end{equation}
As in the case of previous models, the displacement field of this shell combines two maps, namely,
$\mbs{\Sigma}:\mathcal{C}\to \mathbb{R}^3$, governing the displacement of points on the mid-surface,
and $\mbs{\theta}:\mathcal{C}\to \mathbb{R}^2$, describing the fiber infinitesimal rotation. Here, we identify $\mathbb{R}^{2}$ with all vectors in $\mathbb{R}^3$ orthogonal to $\mbs{E}_3$. Based on these two functions, the displacement of the shell can be written as
\begin{equation}
  \label{eq-mindlin-psi}
  \combo{\Sigma}{\theta}
  \defined
  \mbs{\Sigma}(\mbs{\sigma}) + (\mbs{\theta}(\mbs{\sigma})\times \mbs{E}_3)\,\xi \ ,
\end{equation}
a vector field defined on the shell and with a structure that is of the form~\eqref{eq-struct-psi}.
To simplify the expression of the natural inner product of the Mindlin shell, we first note that the gradient of the displacement has the form
\begin{equation}
  \label{eq-mindlin-grad}
  D\combo{\Sigma}{\theta}(\mbs{\sigma},\xi) =
  \pd{\combo{\Sigma}{\theta}(\mbs{\sigma},\xi)}{\sigma_\alpha}\otimes \mbs{E}_{\alpha} +
  \pd{\combo{\Sigma}{\theta}(\mbs{\sigma},\xi)}{\xi}\otimes \mbs{E}_3\ ,
\end{equation}
and we introduce two arbitrary displacement fields for the shell, namely,
\begin{equation}
  \label{eq-mindlin-two-psi}
\combo{\Sigma}{\theta}(\mbs{\sigma},\xi) = \mbs{\Sigma}(\mbs{\sigma}) + (\mbs{\theta}(\mbs{\sigma})\times \mbs{E}_3)\,\xi \ ,
\qquad
\combo{\Upsilon}{\lambda}(\mbs{\sigma},\xi) = \mbs{\Upsilon}(\mbs{\sigma}) + (\mbs{\lambda}(\mbs{\sigma})\times \mbs{E}_3)\,\xi .
\end{equation}
For these fields we obtain
\begin{equation}
  \label{eq-mindlin-inner}
  \begin{split}
      \structinner{\combo{\Sigma}{\theta}}{\combo{\Upsilon}{\lambda}}
    =& \int_{A}
       t
    \left(
       \mbs{\Sigma}(\sigma)\cdot \mbs{\Upsilon}(\sigma) + \ell^2\, \mbs{\Sigma}'(\sigma)\cdot \mbs{\Upsilon}'(\sigma)
       \right) \;\mathrm{d}A
       \\
    &+
\int_{A}
    \left(
      \mbs{\theta}(\sigma)\cdot \mbs{i}_{\mathcal{F}}\,\mbs{\lambda}(\sigma) +
      \ell^2\, \mbs{\theta}'(\sigma)\cdot \mbs{i}_{\mathcal{F}}\, \mbs{\lambda}'(\sigma)
      \right) \;\mathrm{d}A
    \\
    &+
    \int_{A}
      t \,\ell^2\;
      \mbs{\theta}(\sigma)\cdot
      (\mbs{I} + \mbs{E}_3\otimes \mbs{E}_{3})\mbs{\lambda}(\sigma)
      \;\mathrm{d}A
    \ ,
  \end{split}
\end{equation}

As in the case of the Timoshenko beam, the potential energy of an elastic Mindlin shell is
a functional $V_M$ that includes the internal energy due to all its
deformation modes, as well as contributions from the external loading. The equilibrium
displacements $\Sigma, \theta$ are the minimizers
\begin{equation}
  \label{eq:pot_energy_mindlin}
  \combo{\Sigma}{\theta} =
  \arg\inf_{\combo{\Upsilon}{\lambda}} V_M(\mbs{\Upsilon},\mbs{\lambda})\ .
\end{equation}


\section{Coupling dissimilar bodies}%
\label{sec-embedding}
We consider next the main novelty of this article, that is, the formulation of problems consisting
of general \emph{composite solids}, meaning those consisting of a deformable solid medium with
inserted structures that can be one-, two-, or three-dimensional.

The main claim of this article is that many composite bodies can be
analyzed in a unified fashion if they are considered to be continuum bodies with embedded
structural bodies. For this, the inclusions have to be assimilated to abstract structures: fibers to
beams, gravel to rigid bodies, etc. To support this claim we will resort to the Arlequin method, a
procedure introduced to formulate boundary value problems of independent deformable bodies in such a way that their kinematics are made compatible. In the sequel, we will extend this method to the kind of
composite solids described before.

\subsection{The Arlequin method}%
\label{subs-arlequin}
The Arlequin method is a framework proposed in the early 2000s to glue different models
in small strain solid mechanics~\cite{dhia2001iu,dhia2005iu}. The formulation of the method
includes constraints that enforce the compatibility of the displacement fields in the shared
domain but in a way that is dictated by the energy of the solution. In small strain
solid mechanics, where the energy is equivalent to the $H^1$ norm of the solution, the compatibility
is enforced in the $H^1$ sense, using the inner product of this Hilbert space. In the developments that follow, we will continue to enforce the compatibility between the continuum medium and the structure using the corresponding energy inner product, as identified in Section~\ref{sec-model}. For each structural type, these inner products are the ones that are naturally dictated by the geometry and functional spaces of the structural model.

We start by recalling the Arlequin method. For that, consider two elastic bodies occupying (open, bounded) domains
$\Omega_1,\Omega_2\subset \mathbb{R}^3$, respectively, and modeled with small strain kinematics. If
both bodies are subject to boundary conditions of displacement and traction, as well as (possibly)
volume forces, their equilibrium configurations are the displacement fields $\mbs{u}_1\in X_1 \defined
[H^1_\partial(\Omega_1)]^3$ and $\mbs{u}_2\in X_2 \defined [H^1_\partial(\Omega_2)]^3$ that minimize the potential energies of the bodies, respectively. That is,
\begin{equation}
  \mbs{u}_1 = \arg\inf_{\mbs{u}\in X_1}V_{LE,1}(\mbs{u})\ ,
  \qquad
  \mbs{u}_2 = \arg\inf_{\mbs{u}\in X_2}V_{LE,2}(\mbs{u})\ .
\end{equation}
The solution spaces $X_{\alpha}$, $\alpha=1,2$, are the Hilbert spaces of vector fields defined in $\Omega_{\alpha}$ with components in $H^1(\Omega_{\alpha})$ and satisfying the Dirichlet boundary conditions on the corresponding part of the boundary.

The Arlequin method enforces $\mbs{u}_1$ and $\mbs{u}_2$ to be compatible on $\tilde{\Omega} =
\Omega_1\cap\Omega_2\neq\emptyset$ by solving the constrained variational problem
\begin{equation}
  \label{eq-arlequin}
  \inf_{\mbs{u}_1\in X_1, \mbs{u}_2\in X_{2}}
  \sup_{\mbs{\gamma}\in H^1(\tilde{\Omega})}
  \left[
    V_{LE,1}(\mbs{u}_1) + V_{LE,2}(\mbs{u}_2) +
    \langle \mbs{\gamma} , \mbs{u}_1 - \mbs{u}_2 \rangle_{\tilde{\Omega}}
    \right] \ ,
\end{equation}
where $\langle\cdot,\cdot\rangle$ is the $H^1$ inner product defined in Eq.~\eqref{eq-inner-linearized}.
In Eq.~\eqref{eq-arlequin}, $\mbs{\gamma}\in [H^1(\tilde{\Omega})]^3$ plays the role of a Lagrange multiplier,
and we will refer to it as such. Note, however, that in the classical theory of constrained
optimization,  Lagrange multipliers are used to impose restrictions by pairing them, in the sense
of $L^2$, with the constraints. Remarkably, the Arlequin problem is well-posed and its (mixed)
finite element discretization, too (see \cite{dhia2001iu}).

In the original Arlequin method, the contributions of the two tied bodies in the overlap region are
weighed to avoid the energy in $\tilde{\Omega}$ from being included twice
\cite{dhia2001iu,dhia2005iu}. In this work, we assume that the potential
energy is dominated by one of the models in the overlap region which furthermore, is small. There is
a modeling error in these approximations, which we will ignore.

\subsection{Abstract formulation of compatible composite bodies}
\label{subs-abstract}
In this section, we extend the Arlequin method to composite bodies. To do so, we apply the ideas of Section~\ref{subs-arlequin} to a pair of dissimilar bodies: one is a continuum deformable solid and the other one is an abstract structural model.

\begin{figure}[ht]
    \centering
    \includegraphics[width=0.7\textwidth]{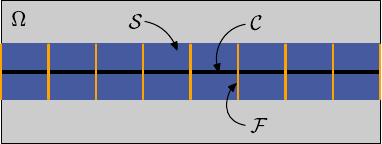}
    \caption{Embedded problem schematic. The elastic body, $\Omega$, is represented in gray, the structure, $\mathcal{S}$, in blue and its base, $\mathcal{C}$, and fiber, $\mathcal{F}$, in black and yellow respectively.}
    \label{fig-embedded-scheme}
\end{figure}

First, let us consider an elastic body as defined in Section~\ref{subs-linear-elas}. Its
displacement field is a function $\mbs{u}:\Omega\to \mathbb{R}^3$ that minimizes the potential
energy $V_{LE}$. Then, let us consider a structure $\mathcal{S}$ with projection
$\Pi:\mathbb{R}^3\to \mathcal{S}$ and generalized displacement field $\combo{\Sigma}{\Theta}$. These
two fields minimize the potential energy $V_{ST}$ in the absence of any interaction with the elastic body. Figure~\ref{fig-embedded-scheme} illustrates the setting.

To enforce the compatibility of the displacements between the continuum body and the structure, we assume that they share some region of the space
\begin{equation}
  \label{eq-btilde}
  \emptyset \neq
  \tilde{\mathcal{S}} \defined \Omega \cap \mathcal{S}\ ,
  \qquad
  \tilde{\mathcal{C}} \defined \Pi_{\mathcal{C}}(\tilde{\mathcal{S}})\ ,
  \qquad
  \tilde{\mathcal{F}} \defined \Pi_{\mathcal{F}}(\tilde{\mathcal{S}})\ .
\end{equation}
Then, generalizing the variational expression~\eqref{eq-arlequin}, we search for the solution
\begin{equation}
  \label{eq-saddle}
  \inf_{\mbs{u}\in X, \combo{\Sigma}{\Theta}\in Y}
\sup_{\mbs{\gamma}\in Z}
\left[
    V_{LE}(\mbs{u}) +
    V_{ST}(\Psi) +
    \structinner{\mbs{\gamma}}{\combo{\Sigma}{\Theta}-\mbs{u}\circ \mathcal{E}}
    \right]
 \ ,
\end{equation}
where the inner product has to be restricted to~$\tilde{\mathcal{S}}$, and the space for the
multipliers, denoted as $Z$, has yet to be defined. Henceforth, and for the sake of notational simplicity, we assume --- without loss of generality --- that the entire base of the structure is  embedded in the solid, i.e. $\tilde{\mathcal{C}} = \mathcal{C}$.


There are several theoretically valid choices for $Z$. The simplest one would be to select
$\gamma=\Gamma\circ \mathcal{E}$, and $\Gamma\in H^1(\tilde{\mathcal{S}})$. However, for reasons
that will become apparent in Section~\ref{subs-wellposed}, we select $Z\equiv Y$. Thus, and
following the notation employed
in Section~\ref{subs-reduced}, a multiplier is a generalized field $\combo{\gamma}{\mu}\in Y$, i.e.,
\begin{equation}
  \label{eq-multiplier}
  \combo{\gamma}{\mu}(\sigma,\xi) =
  \gamma(\sigma) + \mu(\sigma) * \xi\ .
\end{equation}

\subsection{Well-posedness of the linked problem}
\label{subs-wellposed}
In this section, we prove that all the problems of tied continuum/structures abstractly defined in
Section~\ref{subs-abstract} are well-posed. We restrict our analysis to the case where the
structure has no Dirichlet boundary conditions: the existence and uniqueness of its solution will
depend on the way it is linked to the continuum body, which is assumed to have some boundary with
constrained displacements.

To present our analysis, let us write problem~\eqref{eq-saddle} in the
standard form of a \textit{mixed problem}: find the displacements $u\in X,\ \combo{\Sigma}{\Theta}\in Y$ and the
multipliers $\combo{\gamma}{\mu}\in Z$ such that
\begin{equation}
  \begin{split}
  \label{eq-mixprob}
  a(u,\Sigma, \Theta ; v, \Gamma, \beta )
  + b(v, \Gamma,\beta; \gamma, \mu)
  & = f(v, \Gamma, \beta )\ , \\
  b(u, \Sigma, \Theta; \alpha, \zeta)
  & = 0\ .
  \end{split}
\end{equation}
for all $v\in X,\ \combo{\Gamma}{\beta}\in Y$ and $\combo{\alpha}{\zeta}\in Z$. Here, the bilinear forms $a(\cdot , \cdot)$ and the
force term $f(\cdot)$ collect, respectively, the bilinear forms and the
force terms of the continuum solid and the structure. The bilinear form $b(\cdot,\cdot)$ is
responsible for the constraint. Specifically, we employ
\begin{equation}
  \begin{split}
  \label{eq:bilinearF}
  a(u,\Sigma,\Theta; v, \Gamma, \beta)
    &\defined a_{LE}(u , v) + a_{ST}(\Sigma,\Theta ;  \Gamma,\beta)\ ,
      \\
    b( u, \gamma, \mu ; \alpha, \zeta)
    &\defined
      \structinner{\combo{\alpha}{\zeta}}{\combo{\gamma}{\mu} - u\circ \mathcal{E}}\ ,
      \\
  f (v, \Gamma, \beta)
  &\defined f_{LE}(v) + f_{ST}( \Gamma,\beta )\ .
  \end{split}
\end{equation}

To analyze problem~\eqref{eq-mixprob}, we need to define the appropriate norms in all the functional
spaces that we have introduced for the three-dimensional solid and the structure. First, given $\combo{\Sigma}{\Theta}\in Y$, we define
\begin{equation}
  \label{eq-triplenorm}
  \triplenorm{\combo{\Sigma}{\Theta}}
  \defined \structinner{\combo{\Sigma}{\Theta}}{\combo{\Sigma}{\Theta}}^{1/2}\ ,
\end{equation}
and recall the standard norm for elements in $Y$
\begin{equation}
  \label{eq-hilbert-norm}
  \| \combo{\Sigma}{\Theta} \|_Y
  \defined
  \left(
    \| \Sigma \|_{H^1(\tilde{\mathcal{C}})}^2
    + \ell^2\;
    \| \Theta \|_{H^1(\tilde{\mathcal{C}})}^2
  \right)^{1/2}\ .
\end{equation}
The norm on the solution space is defined, for all triplets $(u,\Sigma,\Theta)\in
X\times Y$, as
\begin{equation}
  \label{eq-norm3}
  \| (u,\Sigma,\Theta) \|_{X\times Y}
  \defined
  \left(
    \| u \|_{H^1(\Omega)}^2 + \| \combo{\Sigma}{\Theta} \|^2_Y
    \right)^{1/2}\ .
\end{equation}
Finally, the subspace
\begin{equation}
  \label{eq-kernel}
  \textrm{ker}\, b
  \defined
  \left\{
    u\in X,\ \combo{\Sigma}{\Theta}\in Y\
    \hbox{such that }
    b(u,\Sigma,\Theta; \alpha,\zeta) = 0
    \hbox{ for all } (\alpha,\zeta)\in Z
  \right\}
\end{equation}
will play a key role in the subsequent analysis.

With this formulation, and using the standard
theory of mixed problems \cite{roberts1989vm,brezzi1991tn,auricchio2005tt}, we can state the following,
\begin{theorem}
  \label{th-general-uniq}
  A mechanical system consisting of a continuum body and a structure with compatible kinematics has displacements that are obtained as the unique stationary solutions of the saddle point problem~\eqref{eq-saddle} if the following conditions are met:
  \begin{itemize}
  \item The bilinear form $a(\cdot,\cdot)$ is continuous in $X\times Y$,
  \item The linear form $f(\cdot)$ is continuous in $X\times Y$,
  \item The bilinear form $a(\cdot,\cdot)$ is elliptic in $\ker b$, i.e., there exists a constant $\alpha>0$ such that
  \begin{equation}
  a(u,\Sigma,\Theta ; u, \Sigma, \Theta)
  \ge \alpha\; \| (u, \Sigma, \Theta ) \|^2_{X\times Y}\ .
\end{equation}
for all $(u,\Sigma,\Theta)\in X\times Y$ in $\textrm{ker}\; b$.

  \item The norm $|||\cdot|||$ on $Z$ is uniformly bounded from below by the norm $\|\cdot\|_Y$.
  \end{itemize}

  \begin{proof}
  The first four conditions of the theorem and the \textit{inf-sup} condition are the necessary requirements
  to prove the statement,  according to the standard theory of mixed problems. To prove the latter,
  for each $\combo{\gamma}{\mu}\in Y$,
  \begin{equation}
    \begin{split}
      \sup_{u\in X, \combo{\Sigma}{\Theta}\in Y}
      \frac{b(u,\Sigma,\Theta; \gamma,\mu)}{\| (u,\Sigma,\Theta) \|_{X\times Y}}
      &\ge
      \sup_{ \combo{\Sigma}{\Theta}\in Y}
        \frac{b(0,\Sigma,\Theta; \gamma, \mu )}{\| \combo{\Sigma}{\Theta} \|_Y}
      \\
      &\ge
        \frac{\structinner{\combo{\gamma}{\mu}}{\combo{\gamma}{\mu}}}{\| \combo{\gamma}{\mu} \|_{Y}}
        \\
    &\ge
    C\; \| \combo{\gamma}{\mu}   \|_{Y}\ ,
    \end{split}
  \end{equation}
  for some $C>0$ independent of the solution. This is precisely  the inf-sup condition required for
  well-posedness of the mixed problem.
  \end{proof}
\end{theorem}

\begin{theorem}
\label{th-mech}
A mechanical system consisting of a continuum body and a structure of the type described in Section~\ref{sec-model}
and such that $\tilde{\mathcal{S}} = \tilde{\mathcal{C}} \times \mathcal{F}$, meaning that every fiber
in the structure is either completely embedded or not embedded at all,
verifies the requirements of Theorem~\ref{th-general-uniq}.

\begin{proof}
The continuity of the bilinear form depends on the specific model. For the case of beams and plates
discussed in section~\ref{sec-model}, these proofs can be found in~\cite{Reddy2002},\cite{ciarlet1988ux}. The second requirement is the continuity of the linear form, which depends on the loads of the structural model and is usually assumed to be continuous.

The ellipticity of the bilinear form in $\ker b$ for the structural models presented in
Section~\ref{sec-model} can be proved using arguments similar to the ones found in the analysis of
the Arlequin method~\cite{dhia2001iu} or solid-to-beam variational links\cite{romero2023jg}. First,
we note that the function $a(u,\Sigma,\Theta; u,\Sigma,\Theta)$ is positive
definite. Since $a_{LE}(u,u)$ is positive semidefinite and $a_{ST}(\Sigma,\Theta; \Sigma,\Theta)$ is positive definite,
if $a(u,\Sigma,\Theta;u,\Sigma,\Theta)=0$ then $(\Sigma,\Theta)$ must be $(0,0)$ and $u$ must be an
infinitesimal rigid body motion. However, the only rigid body deformation in $\ker b$ is $u$ identically zero because $(\Sigma,\Theta)= (0,0)$.

To prove the ellipticity of $a(\cdot,\cdot)$ in $\ker b$, we can proceed by proof by contradiction.
Let's suppose that there is no $\alpha > 0$ such that $a(u,\Sigma,\Theta;u,\Sigma,\Theta) \geq
\alpha\, || (u,\Sigma,\Theta) ||^2_{X\times Y}$ for any $(u,\Sigma,\Theta) \in \ker b$. In that case,
there is a sequence $\left\{ (u_n,\Sigma_n,\Theta_n) \right\}_{n \in \mathbb{N}} \in \ker b$ such that
\begin{equation}
   \label{eq:seq_kerb}
   ||(u_n,\Sigma_n,\Theta_n)||_{X \times Y} = \beta,  \quad
   \lim_{n \to \infty} a(u_n,\Sigma_n,\Theta_n; u_n,\Sigma_n,\Theta_n) = 0\ .
\end{equation}
where $\beta >0$. Since $\beta =  ||(u_n,\Sigma_n,\Theta_n)||_{X \times Y} \geq
||(\Sigma_n,\Theta_n)||_{[H^1(\mathcal{C})]^3\times [H^1(\mathcal{C})]^r}$,
by Rellich's theorem there is a subsequence $\left\{ (\Sigma_m,\Theta_m) \right\}_{m \in \mathbb{N}}$
that converges in $[L^2(\mathcal{C})]^3\times[L^2(\mathcal{C})]^r$ to a function
$(\bar{\Sigma},\bar{\Theta})$. Taking into account the second requirement
of~\eqref{eq:seq_kerb}, this must be a Cauchy sequence in the norm
 \begin{equation}
   \label{eq:seq_cauchy_norm}
    |||| u,\Sigma,\Theta ||||^2 = ||\Sigma,\Theta ||_{[L^2(\mathcal{C})]^3\times[L^2(\mathcal{C})]^r}^2
                                  + a (u_m,\Sigma_m,\Theta_m ; u_m,\Sigma_m,\Theta_m)\ .
  \end{equation}
Using Korn's inequality and the ellipticity of $a_{LE}(u,u)$, we can prove this norm is equivalent to $||\cdot||_{X \times Y}$.
Hence, the sequence is Cauchy with respect to $||\cdot||_{X \times Y}$ and since $X \times Y$ is a Hilbert space,
it converges to $(\bar{u},\bar{\Sigma},\bar{\Theta}) \in \ker b$. From \eqref{eq:seq_kerb}, it
follows that
  \begin{equation}
   \label{eq:seq_converge}
     a(\bar{u},\bar{\Sigma},\bar{\Theta} ; \bar{u},\bar{\Sigma},\bar{\Theta}) =
     \lim_{m \to \infty} a(u_m,\Sigma_m,\Theta_m ; u_m,\Sigma_m,\Theta_m) = 0\ .
 \end{equation}
Since $a(\cdot,\cdot)$ is positive definite in $\ker b$, $(\bar{u},\bar{\Sigma},\bar{\Theta}) = (0,0,0)$ but
  \begin{equation}
   0 = || (\bar{u},\bar{\Sigma},\bar{\Theta}) ||_{X \times Y} =
  \lim_{m \to \infty} a(u_m,\Sigma_m,\Theta_m; u_m,\Sigma_m,\Theta_m) = \beta > 0\ .
  \end{equation}
  This contradiction implies the ellipticity of $a(\cdot;\cdot)$ in $\ker b$.

Finally, we can prove the lower bounded requirement in $Z$ by using~\eqref{eq-general-inner}
  \begin{equation*}
    \label{eq-general-norm}
    \begin{split}
      \structinner{\combo{\Sigma}{\Theta}}{\combo{\Sigma}{\Theta}}_{\tilde{\mathcal{S}}}
    =& \int_{\tilde{\mathcal{C}}}
       |\mathcal{F}|
    \left(
       \mbs{\Sigma}(\sigma)\cdot \mbs{\Sigma}(\sigma) + \ell^2\, \mbs{\Sigma}'(\sigma)\cdot \mbs{\Sigma}'(\sigma)
       \right) \; \mathrm{d}\mathcal{C}
       \\
    &+
\int_{\tilde{\mathcal{C}}}
    \left(
      \mbs{\theta}(\sigma)\cdot \mbs{i}_{\mathcal{F}}\,\mbs{\theta}(\sigma) +
      \ell^2\, \mbs{\theta}'(\sigma)\cdot \mbs{i}_{\mathcal{F}}\, \mbs{\theta}'(\sigma)
      \right) \;\mathrm{d}\mathcal{C}
    \\
    &+
    \int_{\tilde{\mathcal{C}}}
    |\mathcal{F}|
      \ell^2\;
      \mbs{\theta}(\sigma)\cdot
      (\mbs{I} + \mbs{E}^i\otimes \mbs{E}^{i})\mbs{\theta}(\sigma)
      \; \mathrm{d}\mathcal{C}
    \ ,
  \end{split}
  \end{equation*}
where $\Theta(\sigma) = \hat{\theta}(\sigma)$. Since $i_{\mathcal{F}}$ and $I+E^i \otimes E^i$ are positive definite tensors, the result follows.
  \end{proof}
\end{theorem}

As noted earlier, the choice of the space of Lagrange multiplier is critical for the
well-posedness of the formulation. The proof of the \textit{inf-sup} condition demonstrates this
clearly, as it requires $Z \subseteq Y$. In particular, the space of multipliers must include
both a translation as well as a rotation part taking care of the transfer, respectively, of forces
and moments from the body to the structure, and \emph{vice versa}. To see this more clearly, consider
problems of embedded structures where the base manifold does not move but the fiber manifold does
--- e. g., a beam subject to pure torsion or a shell under transverse shear. In these situations, if
the rotation multipliers are missing, the structure will not receive any (virtual) work from the
embedding body, underestimating the work required to deform the composite body.


\section{Finite element discretization}\label{sec-fem}
\subsection{Discrete problem}
The finite element discretization of problem~\eqref{eq-mixprob} is relatively standard, with
one peculiarity: two interpolation spaces must be defined; one for the displacements of the solid,
and a second one for the generalized displacements of the structure --- and the Lagrange
multipliers, which coincide with the latter. The two spaces need not be compatible.

To define the finite element spaces, let $\Omega_h$ and $\mathcal{C}_h$ be triangulations of the
body and the base manifold of the structure, respectively. A node set $\mathcal{N}_{\Omega}$ must be selected on the body and
three-dimensional elements (for example, tetrahedra or hexahedra) are introduced with shape functions $N_{a}: \Omega_h\to [0,1],\ a\in \mathcal{N}_{\Omega}$. In turn, the structure must be also partitioned on finite elements
whose geometry is the one of the base manifold (triangle or quads for shells, lines for beams,
etc.). If a node set $\mathcal{N}_{\mathcal{C}}$ collects the nodes at the vertices of the
structural elements, the attendant shape functions $M_{b}: \mathcal{C}_h\to [0,1],\, b\in
\mathcal{N}_{\mathcal{C}}$ may now be defined.

To describe the multipliers that are responsible for the linking, let us define
$\tilde{\mathcal{C}}_h\subset \mathcal{C}_h$ to be the part of the structure that is inside the
three-dimensional body and $\tilde{\mathcal{N}}_{\mathcal{C}} \defined \mathcal{N}_{\mathcal{C}}\cap
\tilde{\mathcal{C}}_h$. Then, the shape functions for the multipliers are the same ones as for the
structure, but restricted to a smaller set, i.e., $M_{c}: \tilde{\mathcal{C}}_h \to [0,1]
,\, c\in\tilde{N}_{\mathcal{C}}$.

The finite element solution of problem~\eqref{eq-mixprob} consists of the fields
$\left(\mbs{u}_h, \mbs{\Sigma}_h, \mbs{\Theta}_h, \mbs{\lambda}_h, \mbs{\mu}_h\right)
\in \mathcal{V}_h = X_h \times Y_h \times Z_h$ of the form
\begin{subequations}
\label{eq-disc-var}
  \begin{align}
    \mbs{u}_h(\mbs{x}) &=
                               \sum_{a\in \mathcal{N}_{\Omega}} N_{a}(\mbs{x})\, \mbs{u}_{a}\ ,
                               \label{eq-disc-var1}
    \\
    \mbs{\Sigma}_h(\mbs{\sigma}) &=
                         \sum_{b\in \mathcal{N}_{\mathcal{C}}} M_{b}(\mbs{\sigma})\, \mbs{\Sigma}_{b}\ ,
                                                        \label{eq-disc-var2}
    \\
    \mbs{\theta}_h(\mbs{\sigma}) &=
                         \sum_{b\in \mathcal{N}_{\mathcal{C}}} M_{b}(\mbs{\sigma})\, \mbs{\theta}_{b}\ ,
                                                        \label{eq-disc-var3}
    \\
    \mbs{\gamma}_h(\mbs{\sigma}) &=
                               \sum_{c\in \tilde{\mathcal{N}}_{\mathcal{C}}} M_{c}(\mbs{\sigma})\, \mbs{\gamma}_{c}\ ,
                                                              \label{eq-disc-var4}
    \\
        \mbs{\mu}_h(\mbs{\sigma}) &=
                               \sum_{c\in \tilde{\mathcal{N}}_{\mathcal{C}}} M_{c}(\mbs{\sigma})\, \mbs{\mu}_{c}\ ,
                                                              \label{eq-disc-var5}
  \end{align}
\end{subequations}
such that
\begin{equation}
  \begin{split}
    \label{eq-mixprob-fem}
a({u}_h,{\Sigma}_h, {\Theta}_h ; {v}_h, {\Gamma}_h, {\beta}_h)
  + b ( v_{h}, \Gamma_h,\beta_h ; \gamma_h, \mu_h)
  & = f ( v_h, \Gamma_h, \beta_h )\ , \\
  b ( u_h, \Sigma_h, \Theta_h; \alpha_h, \zeta_h)
  & = 0\ ,
  \end{split}
\end{equation}
hold for all $v_h\in X_h, \combo{\Gamma_h}{\beta_h}\in Y_h$ and $\combo{\alpha_h}{\zeta_h}\in Z_h$.

The boundary value problem~\eqref{eq-mixprob-fem} has a saddle point structure and might be analyzed
using classical techniques for mixed methods. These are typically non-trivial and demand a careful
examination of the effect of the bilinear forms $a(\cdot;\cdot)$ and $b(\cdot;\cdot)$ on the spaces
$X^h$ and $Y^h$.

In contrast with the Galerkin discretization of elliptic boundary value problems,
the discrete saddle point problems do not inherit, in general, the well-posedness of the continuum
formulation. In the case of the Arlequin method, however, this property follows directly under the
conditions that $X_h \times Y_h$ be a closed subspace of $X \times Y$ and $Z_h \subseteq Y_h$.
When these conditions are satisfied, as in the proposed formulation, the proof of the well-posedness of the discrete
formulation follows exactly the same steps employed in Theorem~\ref{th-mech}. Details are omitted.

\subsection{Numerical evaluation of the integrals}
\label{subs-integrals}
The variational equations~\eqref{eq-mixprob-fem} involve three types of integrals. The bilinear and linear forms of the solid and structure are evaluated performing
approximated integrals on their corresponding domains, namely, the full solid and the structure base
manifold. The quadratures to approximate these terms are standard and need not be addressed here.

The terms that arise from the coupling functional $b(\cdot;\cdot)$ involve integrals on the structure region $\tilde{\mathcal{S}}$
and require a special treatment that we detail next. To simplify the notation, we assume for the
remaining of this section that the structure is fully embedded in the solid domain, i.e., $\mathcal{S}\equiv\tilde{\mathcal{S}}$.

To explain how integrals over $\mathcal{S}$ are approximated, consider an integrable function
$f:\mathcal{S} \rightarrow \mathbb{R}$ and its integral
\begin{equation}
\label{eq:integral_structure}
I \defined
\int_{\mathcal{S}} f(\mbs{x})\; \mathrm{d}V
\approx
\int_{\mathcal{C}}\int_{\mathcal{F}} (f\circ \mathcal{E})
  (\mbs{\sigma},\mbs{\xi})
  \; \mathrm{d} \mathcal{F} \; \mathrm{d} \mathcal{C}\ ,
\end{equation}
where the approximation is due to the fact that we have ignored the Jacobian that should appear when
changing the integral from $\mathcal{S}$ to $\mathcal{C}\times \mathcal{F}$. The motivation for this
simplification was already addressed in Section~\ref{subs-reduced}.

In order to calculate these integrals, we propose to use quadrature rules on the fiber and the base, each with the corresponding dimension. These are standard approximations in finite element
codes and the value of the quadrature points and their weights are readily available.

To implement these rules, at the beginning of all calculations, and only this time, $Q_{\base}$
quadrature points and weights need to be calculated in the base manifold, which we denote as $(\sigma_{q},W_{q}) \in \base
\times \mathbb{R}$. Next, we consider the
fiber at each point of the manifold with base coordinates~$\sigma_{q}$ and we define a quadrature on it. This amounts to selecting
$P_{\fiber}$ points which are denoted $\xi_{q,p}, p=1,\ldots,P_{\mathcal{F}}$ and weights $w_{p}$.
Using this idea, the integral \eqref{eq:integral_structure} can be approximated as
\begin{equation}
\label{eq:jacobian_discrete}
\int_{\mathcal{S}} f(x)\; \dV
\approx \sum_{q\in Q_{\base}} \sum_{p\in P_{\fiber}} (f\circ\mathcal{E})(\sigma_{q},\xi_{q,p})
\, W_{q}\, w_{p}\ .
\end{equation}

In this work, the quadrature points selected on the base manifold for evaluating $b(\cdot,\cdot)$
have been chosen to be the same as those employed for evaluating the bilinear and linear forms of
the structure itself. This choice has the advantage of utilizing the already-created quadrature points of the base from its mesh.
Additionally, as will be seen later, it leverages a certain condition in the structure's mesh
-— that it must be finer than that of the solid -— in order to generate a complete quadrature rule.


\subsection{Mesh size constraints in tied models}
\label{sec-instability}

It is well-known that the Arlequin method might exhibit numerical instabilities when employed to tie
two finite element solids with dissimilar mesh sizes~\cite{dhia2001iu,Guidault2007}. Existing analyses
show that an appropriate choice of the space of Lagrange multipliers suffices to suppress these
undesirable approximation effects. A similar undesirable pathology afflicts Arlequin embeddings of
mixed-dimensionality, like the ones advocated heretofore. However, as we will discuss next, their
remedy is simple: it suffices to select the space of multipliers to be the same as the space of
generalized displacements of the structural model \emph{and} ensure that the mesh size of the
latter is not coarser than the mesh of the embedding solid in the interface region.

The analysis of the finite element approximation of connected structures with the Arlequin method is
more complicated than in the original method since the approximation spaces in the structure are
different than the one of the solid. Nevertheless, some analytical results can be obtained, and we
discuss here the analysis of embedded structures without rotational degrees of freedom.
no rotations.

\begin{figure}[ht]
    \centering
        \includegraphics[width=0.45\textwidth]{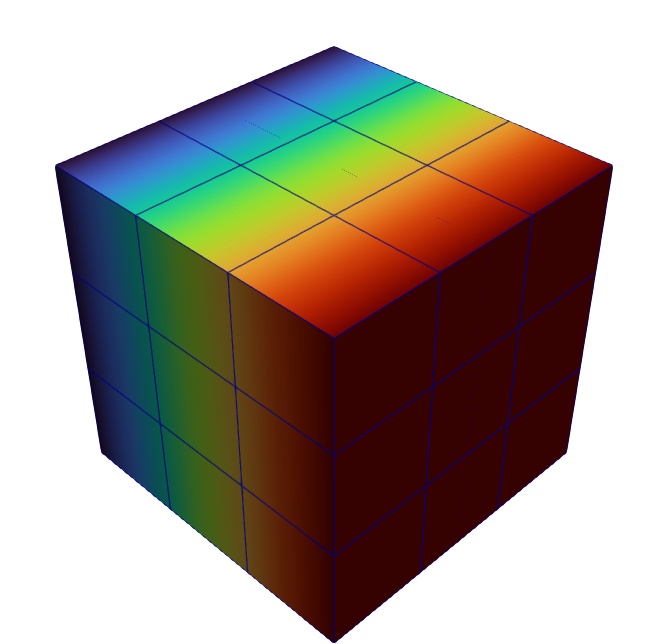}
        \includegraphics[width=0.45\textwidth]{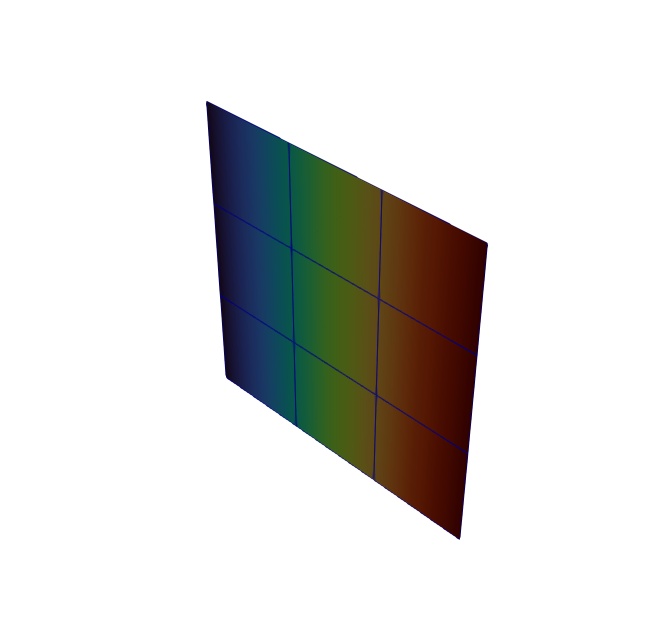}
        \\[1em]
        \includegraphics[width=0.45\textwidth]{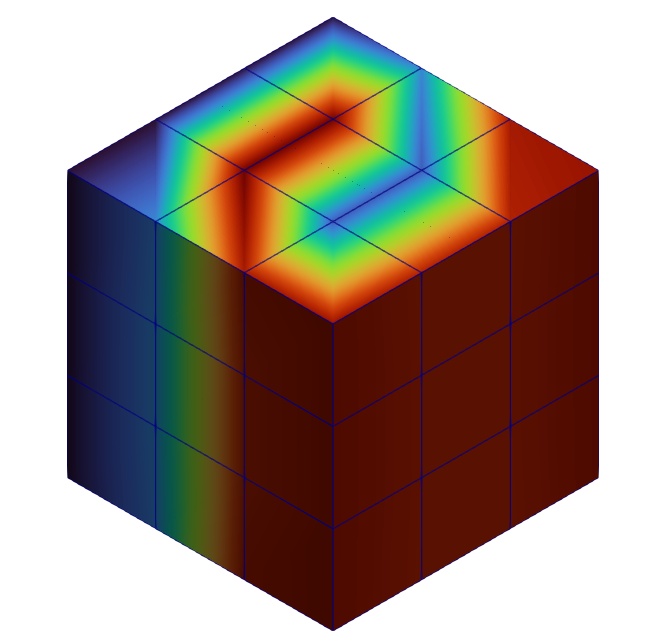}
        \includegraphics[width=0.45\textwidth]{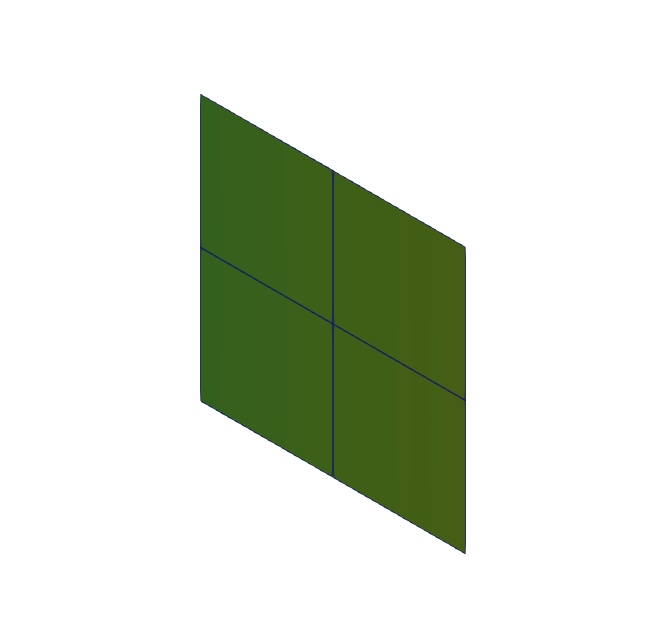}
        \\[1em]
        \includegraphics[width=0.45\textwidth]{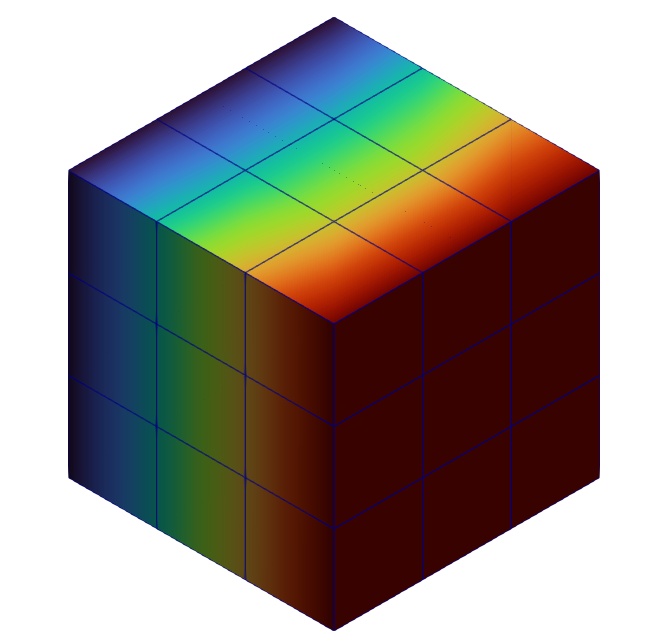}
        \includegraphics[width=0.45\textwidth]{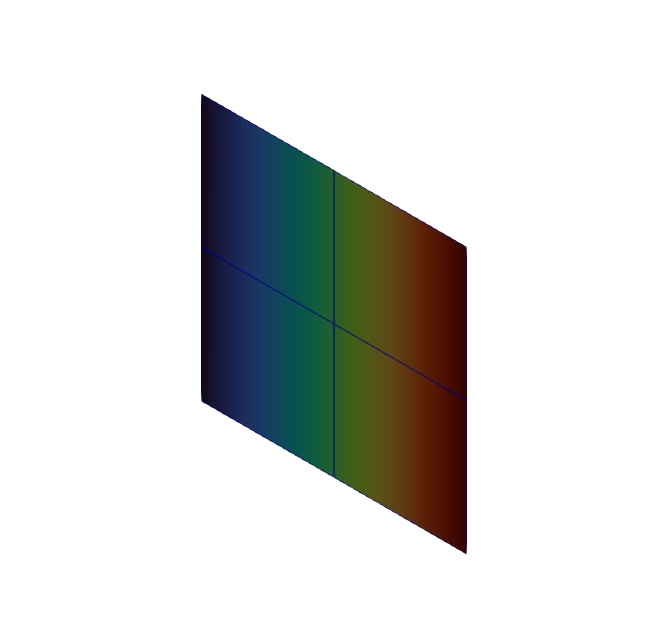}
    \vspace{5mm}
        \includegraphics[width=0.45\textwidth]{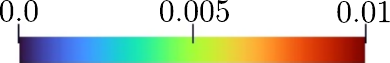}%
    \caption{Displacement field for the traction problem. Above, the problem with the solid and the shell has the same mesh size, $\hs=\hb=\frac{1}{3}$. In the middle, the shell has a larger mesh size $\hb = \frac{1}{2}>\hs=\frac{1}{3}$. Finally, at the bottom, the shell has a larger mesh size, $\hb = \frac{1}{2}>\hs=\frac{1}{3}$, but it is fixed on the left side, like the solid.}
    \label{fig-unstable}
\end{figure}

Before starting the study of the stability issue, let us first denote the mesh sizes
of the solid and the structure as $h_{\Omega}$ and $h_{\mathcal{C}}$, respectively. Next, consider
a solid cube with an embedded shell, as shown in Figure~\ref{fig-unstable}, to showcase the nature
of the instabilities. The cube is subjected to a
displacement on its right face and fixed on the left. As shown in Figure~\ref{fig-unstable} (top),
if the mesh size of the solid and the shell are equal, the correct solution is obtained. In fact,
the same behavior can be observed whenever the mesh size of the shell is smaller than that of the
solid, i.e. $h_{\mathcal{C}} \le h_{\Omega}$. However, if $h_{\mathcal{C}} > h_{\Omega}$ (center of Figure~\ref{fig-unstable}), the solution obtained is physically incorrect.

This instability is a numerical problem that arises from discretization and it depends on the ratio
of solid and shell stiffness, the mesh size, and the order of the interpolating spaces.
Analytically, this problem can be traced down to the value of the coercivity constant and, more specifically, to the lack of coercivity of the structure (which is assumed to
have no Dirichlet boundary conditions) that needs to be compensated with a large enough coercivity
constant of the solid part. In fact, as shown in Figure~\ref{fig-unstable} (bottom), if Dirichlet
boundary conditions are imposed on the plate, the problem identified disappears, irrespective of the
mesh size ratio.

Let us next examine the problem for simple structures that only have translational displacements.
Despite its simplicity, this example is sufficient to demonstrate the numerical instability. To keep
the discussion brief, the main results are presented here, but the details are left for Appendix~\ref{app-instabilities}.

In the continuum problem, the structure bilinear form reduces to
\begin{equation}
a_{ST}(\Sigma_{\hb},\Gamma_{\hb}) = \int_{\mathcal{C}} D \Sigma_{\hb} \cdot C_{\Sigma} \; D\Gamma_{\hb} \textrm{d}\mathcal{C}\ ,
\end{equation}
and the coupling term between solid and structure is
\begin{equation}
\label{eq:bb_trans}
b(\lambda_{\hb},\Gamma_{\hb}-v_h)= \int_{S} \left[ \lambda_{\hb} \cdot (\Gamma_{\hb}-v_h) + \ell^2 D\lambda_{\hb} \cdot D(\Gamma_{\hb} - v_{h} )\right] \dV\ .
\end{equation}
We can then estimate the coercivity constant by identifying $\ker b$ with $u |_{\tilde{S}}=
\Gamma_{\hb}$, hence
\begin{equation}
\label{eq:coercivity_cont_trans}
\begin{array}{lcl}
a(u,\Sigma,\Theta ; u,\Sigma,\Theta)  & \ge &
 \min \left(\epsilon_1, \epsilon_2 \frac{(1-\epsilon_1)}{C_{P_{\Omega}}}\right) \underline{C_{\Omega}} ||u||_{H^1(\Omega)}
 									      \\
 									    &   &
                                           + \min \left( (1-\epsilon_2)\frac{(1-\epsilon_1)\underline{C_{\Omega}}}{C_{P_{\Omega}}} \alpha_{\Sigma}, \underline{C_{\Sigma}}\right) ||\Sigma||_{H^1(\mathcal{C})}^2\ ,
\end{array}
\end{equation}
where $\underline{C_{\Sigma}}$ is a lower bound for the norm of the elasticity tensor,
$C_{P_{\Omega}}$ is Poincare's constant, and $\alpha_{\Sigma}$, $\alpha_{\theta} > 0$ are constants
that depend on the fiber's geometry. The remaining constants, $0<\epsilon_1,\epsilon_2<1$, appear in
the proof after using Young's inequality and allow to compensate the lack of coercivity in the
structure with the excess of the latter in the solid. This reasoning, however, does not extend
directly to the discrete case. Although we can prove ellipticity of the bilinear form
(Theorem~\ref{th-mech}), the calculation of the lower bound constant is not straightforward. In
fact, it depends on the structural model and the interpolation fields we use for solving the
discrete models.

We distinguish here two cases depending on the mesh size of the structure with respect to the solid
and address them separately.

\paragraph{Mesh on the structure is finer than mesh on the solid.} When the structure only has translation degrees of freedom and its mesh is finer than the one in the
solid, we have $\tilde{\Omega}_h \subseteq \tilde{C}_h$, where $\tilde{\Omega}_h$ is the part of the
discrete solid which coincides with the discrete structure. Then, $u_h|_{\tilde{\Omega}_h} \in
\tilde{C}_h$. Taking this into account,  an expression analogous to
Eq.~\eqref{eq:coercivity_cont_trans} can be obtained. As it was noted for the continuum case, the
lack of coercivity in the structure is compensated by the solid's.

\paragraph{Mesh on the structure is coarser then mesh on the solid.} When the mesh on the structure is coarser then the mesh on the solid, the analysis can not follow
the same argument as in the previous case. Now, $\ker b$ is no longer identified with solutions that
are equal in $\tilde{\mathcal{S}}$. Instead, the kernel is identified with the lifting of the
structure being equal to the projection of the solid displacements onto a space with the mesh size
of the structure. Designating as $\uhb$ the projection of the solid displacement $u_h$ onto the functional space given by $\hb$, it can be seen that in $\ker b_{\hb}$, $\Sigma_{\hb} = \uhb$ and we can get the following estimate of the coercivity constant,
\begin{equation}
\label{eq:coercivity_disc_trans}
\begin{array}{lcl}
a(u_{\hs},\Sigma_{\hb};u_{\hs},\Sigma_{\hb})
                                            & \ge & \underline{C_{\Sigma}}
                                            \left(
                                            \min \left(\frac{\epsilon_1 \underline{C_{\Omega}}}{\underline{C_{\Sigma}}}, \epsilon_2  \frac{(1-\epsilon_1) \underline{C_{\Omega}}}{C_{P_{\Omega_{\hs} \; \underline{C_{\Sigma}}}}} \right) ||u_{\hs}||_{H^1(\Omega_h)}^2  \right.
                                            \\
                                           & &  -(1-\epsilon_2)  \frac{(1-\epsilon_1) \underline{C_{\Omega}}}{C_{P_{\Omega_{\hs}}} \; \underline{C_{\Sigma}}} \beta \hb^{2(p+1)}|u_{\hs}|_{H^{p+1}(\Omega_{\hs})}^2 \\
                                         & & \left.
                                         + \min \left( (1-\epsilon_2)  \frac{(1-\epsilon_1) \; |\fiber| \; \underline{C_{\Omega}}}{C_{P_{\Omega_{\hs}}} \; \underline{C_{\Sigma}}}, 1 \right)  ||\Sigma_{\hb}||_{H^1(\mathcal{C}_{\hb})}^2
                                          \right)\ ,
\end{array}
\end{equation}
for $0 < \epsilon_1, \epsilon_2 < 1$, and $\beta$ is a constant that does not depend on the mesh. Finally, $p$ is the order of the finite element approximation of the structure.

There is an important difference between this result and the one for the finer structure mesh. In
the current situation, the coercivity is again subtracted from the solid though the product
$(1-\epsilon_1)(1-\epsilon_2)$, but now its value limited by the error estimate, since that part
limits the coercivity of the solid itself. In this case, the approximation on the structure’s mesh
limits the coercivity, potentially leading to numerical instability depending on the structure-solid
mesh size ratio, the ratio of stiffness in the two bodies, and also on the order of the finite element approximation.

We can conclude that when the mesh size of the structure is greater than that of the solid,
numerical instabilities may arise as a consequence of the coercivity limitation. These issues occur
especially when the structure is stiffer than the solid. Limiting our discretizations to ones where
the structure has a smaller or equal mesh size than the solid avoids the identified instabilities altogether.


\section{Numerical examples}\label{sec-examples}

In this section, we present various tests to demonstrate the capabilities of the proposed framework for structures embedded within solids. Results are provided for both beams and shells.

\subsection{Bending beam problem}
\label{subs-bending}

This first example studies a beam-like body under loading that causes it to bend. It is a simple
benchmark which allows us to compare the solutions obtained with solid elements and embedded bodies,
and also with those in the literature~\cite{steinbrecher2020cm}.

\begin{figure}[ht]
    \centering
    \includegraphics[width=0.5\textwidth]{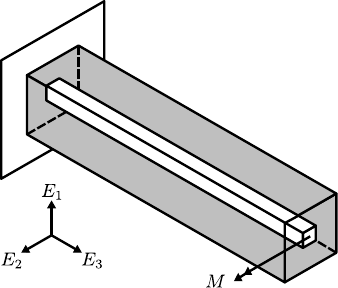}
    \caption{Bending beam problem schematic.}
    \label{fig-bending-beam}
\end{figure}

For the proposed example, we consider a
prismatic body with square base of dimensions $1\times1$ placed on the $(x_1,x_2)$ plane and
length~$5$ along the $x_3$ direction. All the edges are parallel to the coordinate axes. On the
face $x_3=0$, all the body displacements are constrained. On the opposite face, at
$x_3=5$, a distributed surface loading of the form
\begin{equation}
\mbs{t}(x_1,x_2,x_3) = - 12 M x_1 \; \mbs{E_3}\ ,
\end{equation}
is imposed, with $M=-0.0025$. The material of this body is elastic with Young's modulus
$E_{\Omega}=10$ and Poisson's ratio $\nu=0$. The prismatic region in the core of the solid (with square
cross section of side equal to $0.25$) is of a different elastic material, also with zero Poisson's
ratio but with a Young modulus $E_{\mathcal{S}}=5120$. See Figure~\ref{fig-bending-beam} for an
illustration of the composite body.

A first solution to this problem is obtained using a finite element model consisting only of solid
elements, using $16 \times 16 \times 640$ hexahedral elements, for a total of $554,880$ degrees of
freedom. In this approximation, the maximum displacement in the $x_1$ direction has a modulus of
$0.191$ and takes place at the vertices of the face $x_3=5$.

\begin{figure}[ht]
    \centering
    \includegraphics[width=0.45\linewidth]{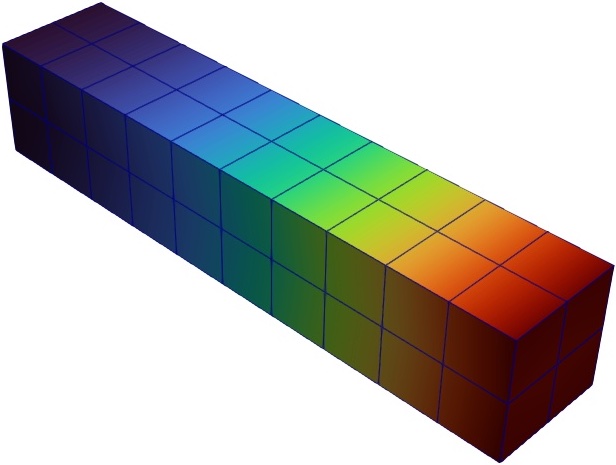}
    \includegraphics[width=0.45\linewidth]{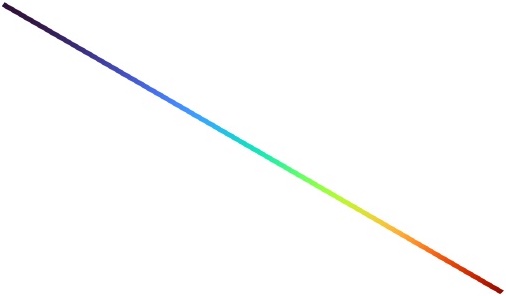}

   \vspace{5mm}

   \includegraphics[width=0.4\textwidth]{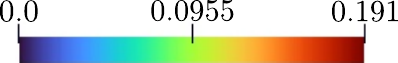}

    \caption{Solid (left) and beam (right) displacement in $x_2$ direction for the bending beam problem.}
    \label{fig-beamBending-sol}
\end{figure}

\begin{figure}[ht]
    \centering
    \includegraphics[width=0.9\textwidth]{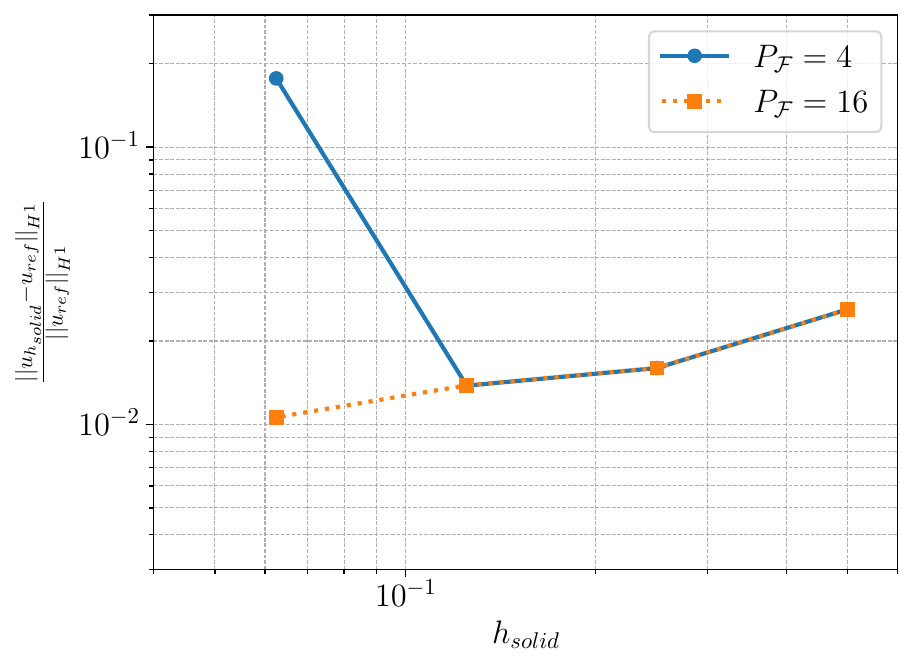}
    \caption{$H^1$ mismatch in the displacement vs. mesh size along the beam's direction for the
      bending problem. Results when
      using $4$ and $16$ quadrature rules for each beam cross section.}
    \label{fig-bending-beam-conv}
\end{figure}

A second solution is obtained by using a single, homogeneous prismatic body of cross
section $1\times1$, with Young's modulus $E_{\Omega}$, and embedding inside it a beam
superposed to the $x_3$ axis, and with Young's modulus $E_{\mathcal{S}}$, cross section area $A=0.25^2$
and inertia $I=0.25^4/12$. Note that this embedded beam is selected to replace the mechanical
effects of the stiffer core of the original body, but it is embedded into a homogeneous body with no
holes. The solid is discretized with a $2\times2\times10$ mesh of hexahedral elements and the beam with $10$ linear elements, leading to a discrete model with 402
unknowns. For this model, the maximum nodal $x_1\textrm{-displacement}$ is 0.185, a merely 3.04\%\ mismatch
relative to first, continuous, model with just 1\%\ degrees of freedom.
Figure~\ref{fig-bending-beam-conv} shows that the mismatch in the $H^1$ norm in the embedded beam
solution is $2.81\%$, relative also to the first model, where the mesh size of the beam is always set to be equal to the one of host solid in the beam's direction. See also Figure~\ref{fig-beamBending-sol} for an illustration of the displacement in the $x_1$
direction for the embedding and embedded bodies in this analysis.

Figure~\ref{fig-bending-beam-conv} shows the errors of the solutions in the two types of
discretizations, as measured in the $H^1$ norm, for various solid mesh sizes, and two different
choices in the number of quadrature points on each fiber. It is important to note that in the case
with fewer quadrature points along the fiber, the solution deteriorates when the mesh size of
the solid becomes comparable to the size of the fiber. This is because the quadrature points must be
able to capture all the solid elements lying within the core region in order to accurately integrate
the response therein. Otherwise, the elements that are not accounted for result in a model with different
stiffness and the outcome is affected. As shown in this figure, if more quadrature points are added along the fiber so that all solid elements are considered during integration (which does not increase the number of degrees of freedom), the expected result is recovered.

The two composite bodies --- the inhomogeneous solid and the one with the embedded beam ---
are not identical, and the same displacement should not be
expected from both. The second one is much easier to discretize, much cheaper in terms of unknowns,
but ignores the overlap of the two deformable bodies on the embedded region. Additionally, the kinematics in the
embedded beam is constrained by Timoshenko's \emph{ansatz}. The errors in
Figure~\ref{fig-beamBending-sol} indicate that both solutions become close when the discretizations
become finer, but one should not expect the convergence of one to another as $h\to0$.

\subsection{Torsion beam problem}
\label{subs-torsion}

Similar to the previous one, this second problem studies the torsion of a prismatic body and can be
compared with results available in the literature \cite{steinbrecher2020cm, khristenko2021cma}. The
geometry, materials, and Dirichlet boundary conditions are the same as in the example of
Section~\ref{subs-bending}. The loading in this second example consist of surface tractions on the face
$x_3=5$ with a distribution given by
\begin{equation}
\tau = 0.05 \left( -x_2 \; E_{1} + x_1 \; E_{2}\right)\ ,
\end{equation}
resulting in a pure torsion moment of modulus equal to $M=1.66 \cdot 10^{-2}$.

A reference solution is obtained using $40\times40\times80$ solid hexahedral elements. On the other hand, a second solution is computed for a prismatic solid with an embedded beam, as in Section~\ref{subs-bending}. In this latter case, the solid is discretized with $4\times4\times8$ hexahedral elements, while the beam is modeled with $8$ linear elements.

This problem is particularly interesting for the method developed in this work since the
deformation of the body leaves the axis of the centroid unaltered, while the fibers --- the cross-sections ---
undergo a rotation around this line. It illustrates a situation in which the translation
Lagrange multipliers are inactive, and the coupling between the solid and the beam is due solely to the rotational Lagrange multipliers.

\begin{figure}[ht]
    \centering
    \includegraphics[width=0.45\linewidth]{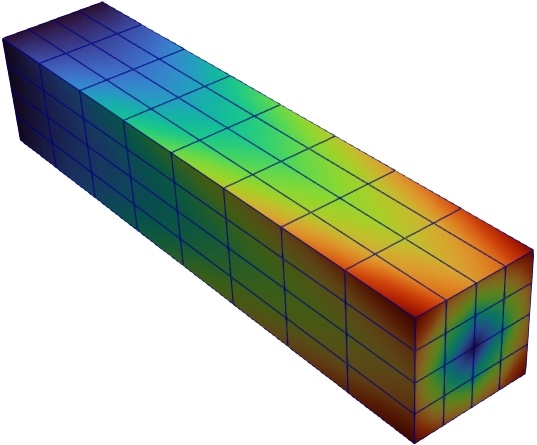}
    \includegraphics[width=0.45\linewidth]{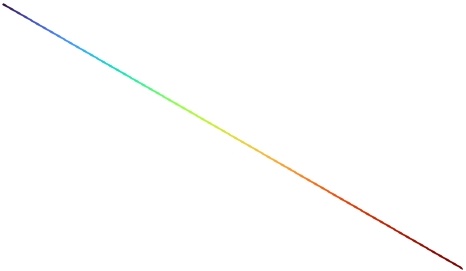}
   \vspace{5mm}
   \includegraphics[width=0.45\linewidth]{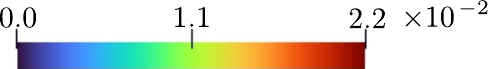}
        \includegraphics[width=0.45\linewidth]{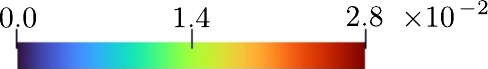}
    \caption{Solid displacement magnitude (left) and beam torsion rotation (right) fields for embedded torsion bending problem.}
    \label{fig-beamTorsion-sol}
\end{figure}

Figure~\ref{fig-beamTorsion-sol} shows the torsion solution, where the displacement magnitude and
the beam’s rotation in the torsion direction can be observed. Note that the beam does not translate
and only rotates, as indicated, since the rotational Lagrange multiplier transfer the required load from the solid. Without this Lagrange multiplier, there would be no compatibility of motions, and the solid would move without accounting for the beam.

\begin{figure}[ht]
    \centering
    \includegraphics[width=0.9\textwidth]{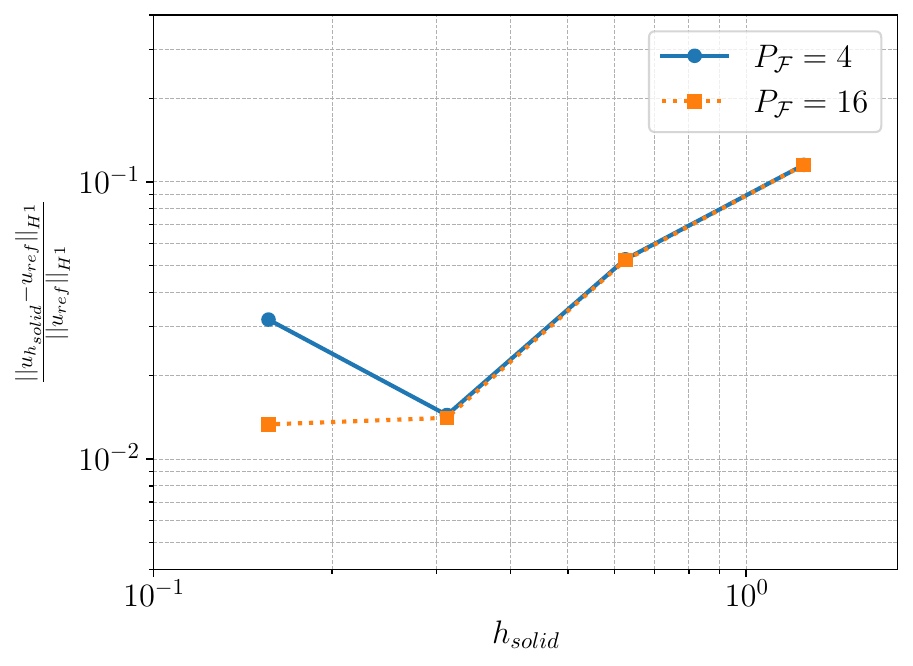}
    \caption{$H^1$ displacement mismatch vs. mesh size along the beam's direction for the torsion beam problem. Results when
      using $4$ and $16$ quadrature rules for each beam cross section.}
    \label{fig-torsion-beam-conv}
\end{figure}

Figure~\ref{fig-torsion-beam-conv} shows an error analysis using the discretization with
$40 \times 40 \times 80$ hexahedral elements as a reference solution. As in the previous case, when
only a few quadrature points are used on the fiber, the solution becomes excessively flexible once
the solid mesh size reaches the core section thickness.

\begin{figure}[ht]
    \centering
    \includegraphics[width=0.5\textwidth]{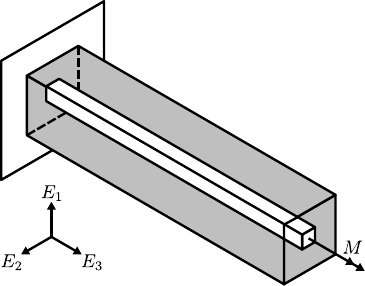}
    \caption{Torsion beam problem schematic.}
    \label{fig-torsion-beam}
\end{figure}

\subsection{Bending shell problem}
\label{subs-shell-bending}

\begin{figure}[ht]
    \centering
    \includegraphics[width=0.5\textwidth]{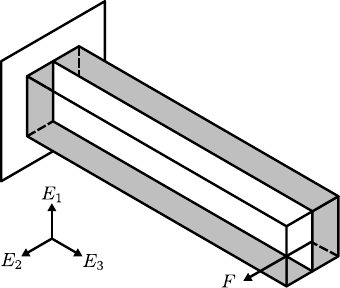}
    \caption{Bending shell problem schematic.}
    \label{fig-bending-shell}
\end{figure}

The first example involving an embedded shell consists of the analysis of a body under loads that
force it to bend. The embedding solid is the same as in the examples of Sections~\ref{subs-bending}
and~\ref{subs-torsion}; the embedding now consists of a shell of thickness $0.25$, parallel to the
$x_1-x_3$ plane and containing the axis of the prism. The Young modulus of the solid and the shell are, respectively,
$E_{\Omega}=10$ and $E_{\mathcal{S}}=640$ and the Poisson's ratio of both materials is zero. A
constant, distributed surface load is applied on the free face of the solid, opposite to the support, and on the $x_3$ direction, with a total resulting force equal to $F=5\cdot 10^{-3}$. Figure~\ref{fig-bending-shell} depicts a schematic of the problem.

A reference solution is obtained using $16\times16\times640$ solid hexahedral elements. On the other hand, a second solution is computed for a prismatic solid with an embedded shell of thickness $0.25$. In this latter case, the solid is discretized with $4\times4\times20$ hexahedral elements, while the shell is modeled with $4\times20$ bilinear quadrilateral elements.

\begin{figure}[ht]
    \centering
    \includegraphics[width=0.45\linewidth]{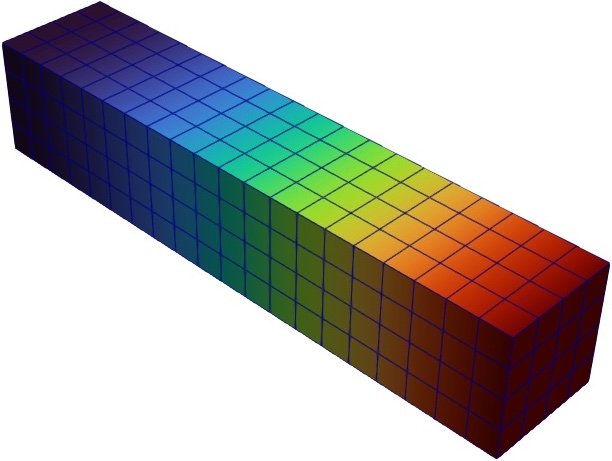}
    \includegraphics[width=0.45\linewidth]{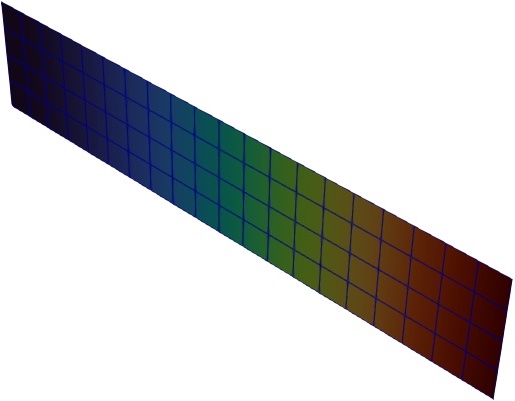}

   \vspace{5mm}
    \includegraphics[width=0.45\textwidth]{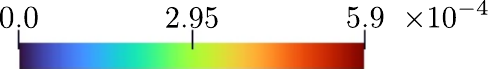}

    \caption{Solid (left) and shell (right) displacement field on the $x_2$ direction for the embedded shell bending problem.}
    \label{fig-shellBengind-sol}
\end{figure}

In Figure~\ref{fig-shellBengind-sol}, the displacement field on the $x_2$ direction is shown for
both the solid and the shell. The compatibility of the displacement between
solid and shell is apparent. The maximum displacement in $x_2$ direction is $0.124$ at the vertices of the face $x_3=5$, corresponding to a $2.78\%$ error relative to the pure solid solution.

\begin{figure}[ht]
    \centering
    \includegraphics[width=0.9\textwidth]{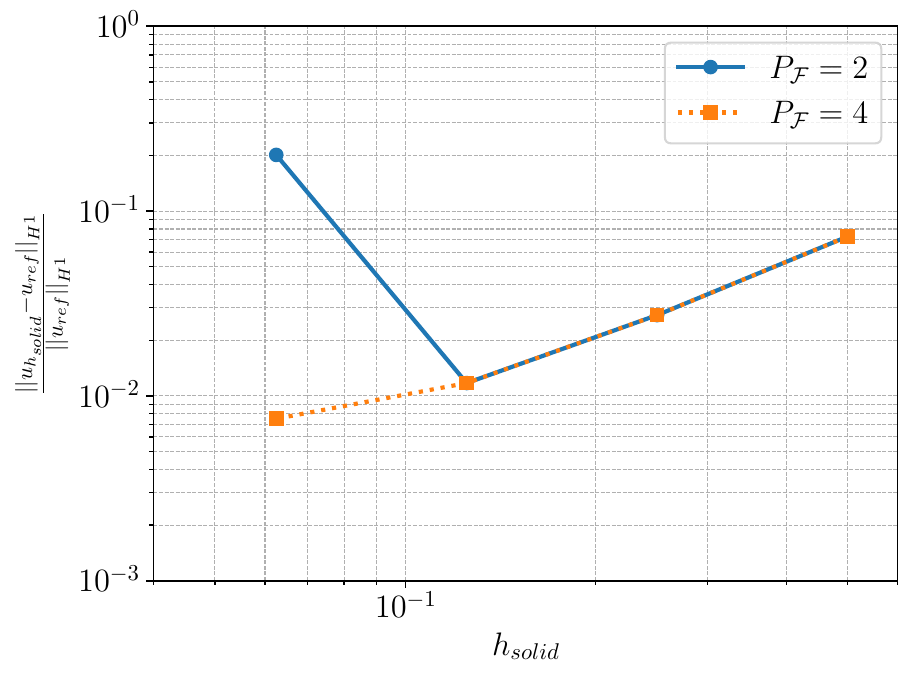}
    \caption{$H^1$ displacement mismatch versus solid mesh size for the bending shell problem. Results when
      using $2$ and $4$ quadrature rules for each shell fiber.}
    \label{fig-bending-shell-conv}
\end{figure}

As in the examples with embedded beams, a mesh sensitivity analysis has also been carried out using a solution
obtained only with solid elements as reference. Figure~\ref{fig-bending-shell-conv} shows this
analysis for two different configurations of quadrature points along the fibers: $2$ and $4$. The
same phenomenon observed in the beam problem can be seen here. When the solid mesh is fine enough,
two quadrature points along the fiber are insufficient to account for all the solid elements falling
within the shell region. When this happens, the solution becomes more flexible (in this case, the
shell is stiffer than the solid), and the results diverge. By increasing the number of quadrature
points on the shell fiber so that all solid elements in the shell region are accounted for, the
solution with the embedded shell becomes closer to the reference. As in the examples of Section~\ref{subs-bending} and \ref{subs-torsion}, we insist
that we do not expect the embedding solution to be identical to the one obtained only discretizing
the continuous solid, given the simplifying assumptions made for the shell and the embedding itself.

\subsection{Shear shell problem}
\label{subs-shell-shear}
Next, we consider a block with an embedded shell under shear loading. The system is designed so
that the deformation of the ensemble does not modify the midsurface of the shell, only inducing the
rotations of its fibers.

\begin{figure}[ht]
    \centering
    \includegraphics[width=0.6\textwidth]{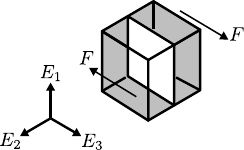}
    \caption{Shear shell problem schematic.}
    \label{fig-shear-shell-scheme}
\end{figure}

The setup consists of a solid cube of size $1 \times 1 \times 1$, with edges parallel to the
Cartesian axes, and a flat shell, parallel to the $(x_1,x_3)$ plane, passing through the center of the cube, and with thickness equal to $0.125$. Both
the solid and the shell are elastic isotropic. The Young's modulus of the solid is $E_{\Omega} =
10$, while the shell material is four times stiffer, with $E_{\mathcal{S}} = 40$. Both bodies have
zero Poisson’s ratio. See Figure~\ref{fig-shear-shell-scheme} for an illustration.

Shear surface loads are applied on the faces $x_2=0$ and $x_2=1$ in opposite direction, with expression
\begin{equation}
\tau |_{x_2=0} = 0.01 \; E_3, \quad \tau |_{x_2=1} = -0.01 \; E_3\ ,
\end{equation}
so a total force $F=0.01$ is applied on each face. Finally, the displacements in directions $x_1$
and $x_2$ are constrained. Since the problem defined in this way is well-posed up to a rigid body motion, the outer nodes at $x_2 = 0.5$ are also fixed. Figure~\ref{fig-shear-shell-scheme} shows a schematic of the problem.

A reference solution is obtained using $64\times64\times64$ solid hexahedral elements. On the other hand, a second solution is computed for a prismatic solid with an embedded beam, as in Section~\ref{subs-bending}. In this latter case, the solid is discretized with $16\times16\times16$ hexahedral elements, while the shell is modeled with $16\times16$ bilinear quadrilateral elements.

As in the case of torsion in beams, this is a problem where the base manifold of the structure ---
in this case, its midsurface --- does not move under the applied loading, yet the fibers rotate. The
shear deformation of the solid is transferred to the shell by means of the rotational Lagrange
multiplier, since the translational one does not play any role.

\begin{figure}[ht]
    \centering
    \includegraphics[width=0.6\textwidth]{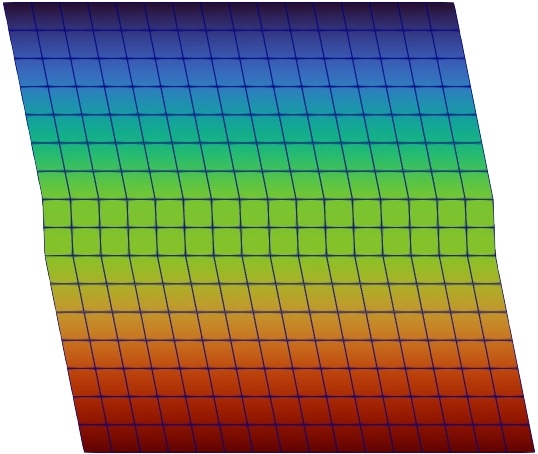}
     \vspace{5mm}
    \includegraphics[width=0.45\textwidth]{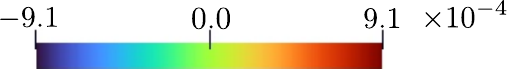}
    \caption{Deformed shape of the shear shell problem, with displacements scaled by a factor of
      100.}
    \label{fig-shear-shell-u}
\end{figure}

Figure~\ref{fig-shear-shell-u} shows the deformed shape of the body and shell, with displacements
scaled by a factor of $100$. The region occupied by the shell, which is four
times stiffer, deforms significantly less than the rest of the much more flexible solid. The only
motion occurring in the shell is that of the director, which has rotated $3.7\times 10^{-4}$ around the $x_2$
direction.

\begin{figure}[ht]
    \centering
    \includegraphics[width=0.9\textwidth]{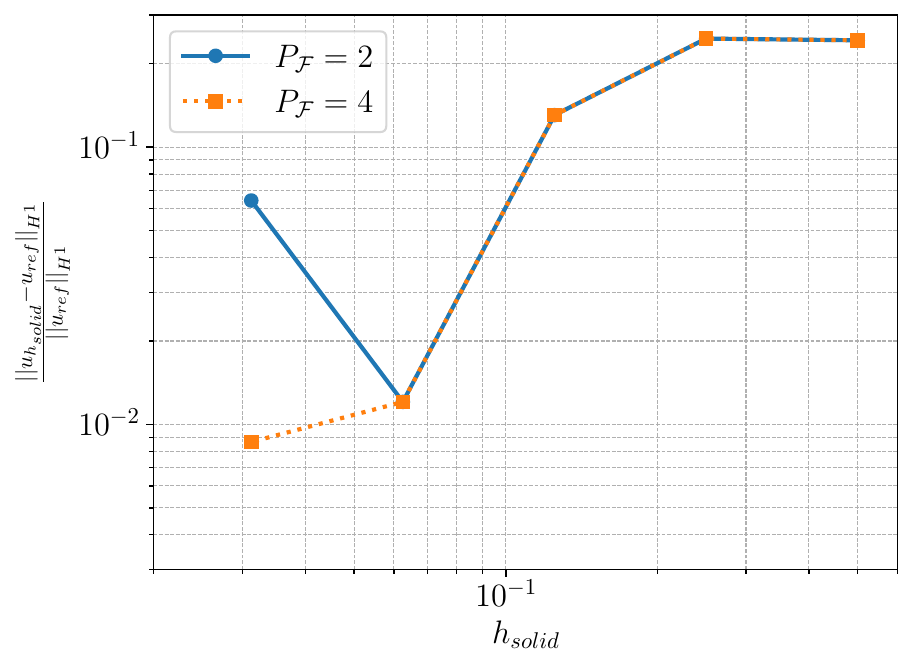}
    \caption{$H^1$ displacement mismatch versus solid mesh size for the shear shell problem. Results when
      using $2$ and $4$ quadrature rules for each shell fiber.}
    \label{fig-shear-shell-conv}
\end{figure}

Figure~\ref{fig-shear-shell-conv} illustrates the mesh sensitivity analysis comparing the solution
obtained with the embedded shell and a reference one obtained only withe solid elements but a much
finer mesh of $64 \times 64 \times 64$ hexahedra. As in previous examples, the effect of an inexact integration is also evident when the solid mesh size reaches that of the fiber. Moreover, it is relevant in this case that the error does not drop below $10\%$ until the mesh size reaches the order of the shell thickness. This is due to the fact that the structure becomes overly stiff with larger mesh sizes, as it affects nodes that are not contained within the volume of the shell. In a way, it is the opposite effect to that observed in all the cases involving incomplete quadrature in the fibers.

\subsection{Elastic solid with reinforcing fibers}
\label{subs-reinforced}
In this last example, a proof of concept of one of the most salient capabilities of the presented
framework is demonstrated. In it, we calculate the effective Young’s modulus of an elastic material reinforced
with fibers, opening an alternative route to standard methods in computational micromechanics.

For this purpose, a unit cubic cell with sides of length~1 is discretized with $10\times10\times10$
hexaedral elements. Periodic boundary conditions are imposed on all the faces of the cell
and imposed displacement of $\pm 0.005$ on faces at $x_3 = \pm 0.5$, respectively (i.e. $\bar{\varepsilon}_{33} = 0.01$). The solid consists of an isotropic
elastic material with Young’s modulus, $E_{\Omega}=1$, and zero Poisson’s ratio, and the fibers are
beams with a length of $0.2$ and a circular cross section of radius $0.025$, a Young’s modulus ten
times larger than the solid, $E_{\mathcal{S}}=10$, and also Poisson’s ratio equal to zero.

\begin{figure}[ht]
    \centering
    \includegraphics[width=0.45\linewidth]{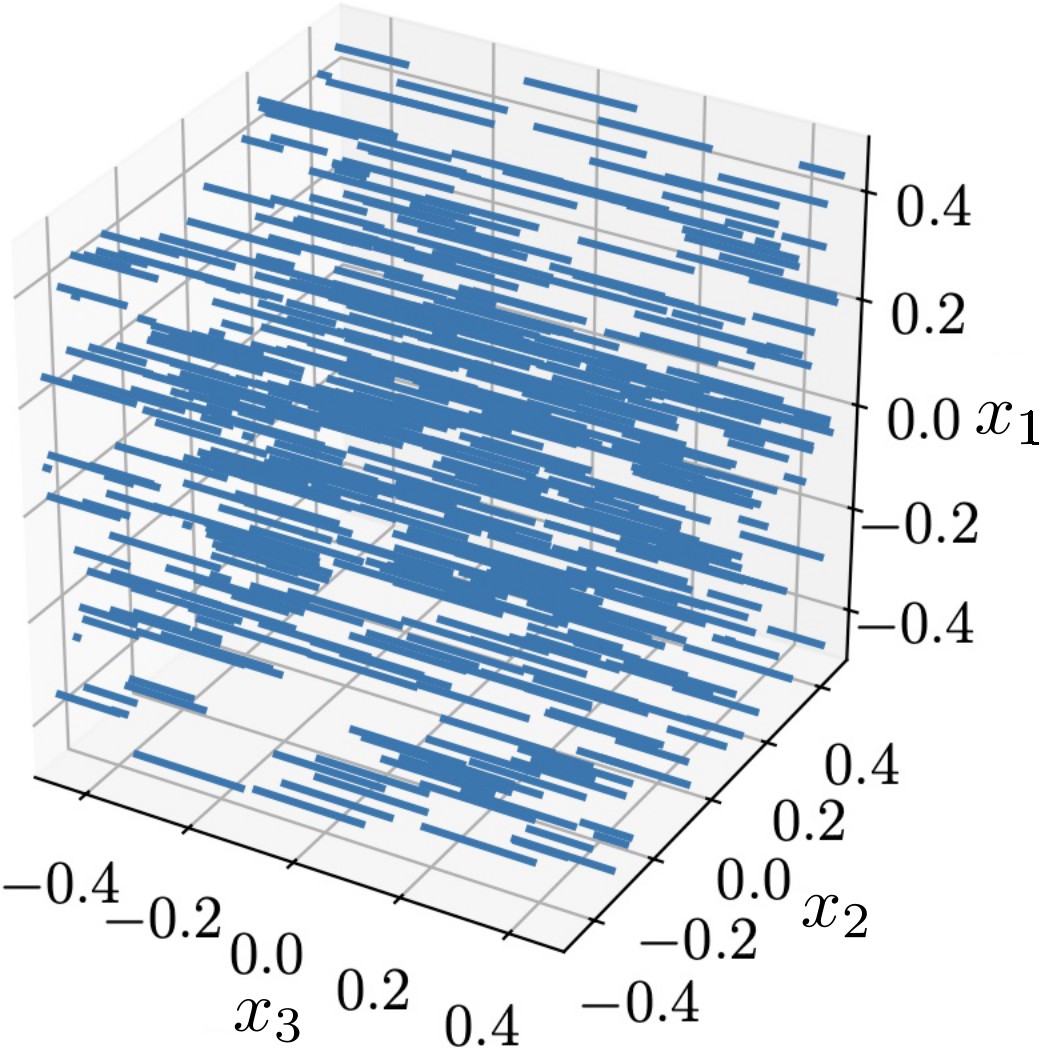}\hfill
    \includegraphics[width=0.45\linewidth]{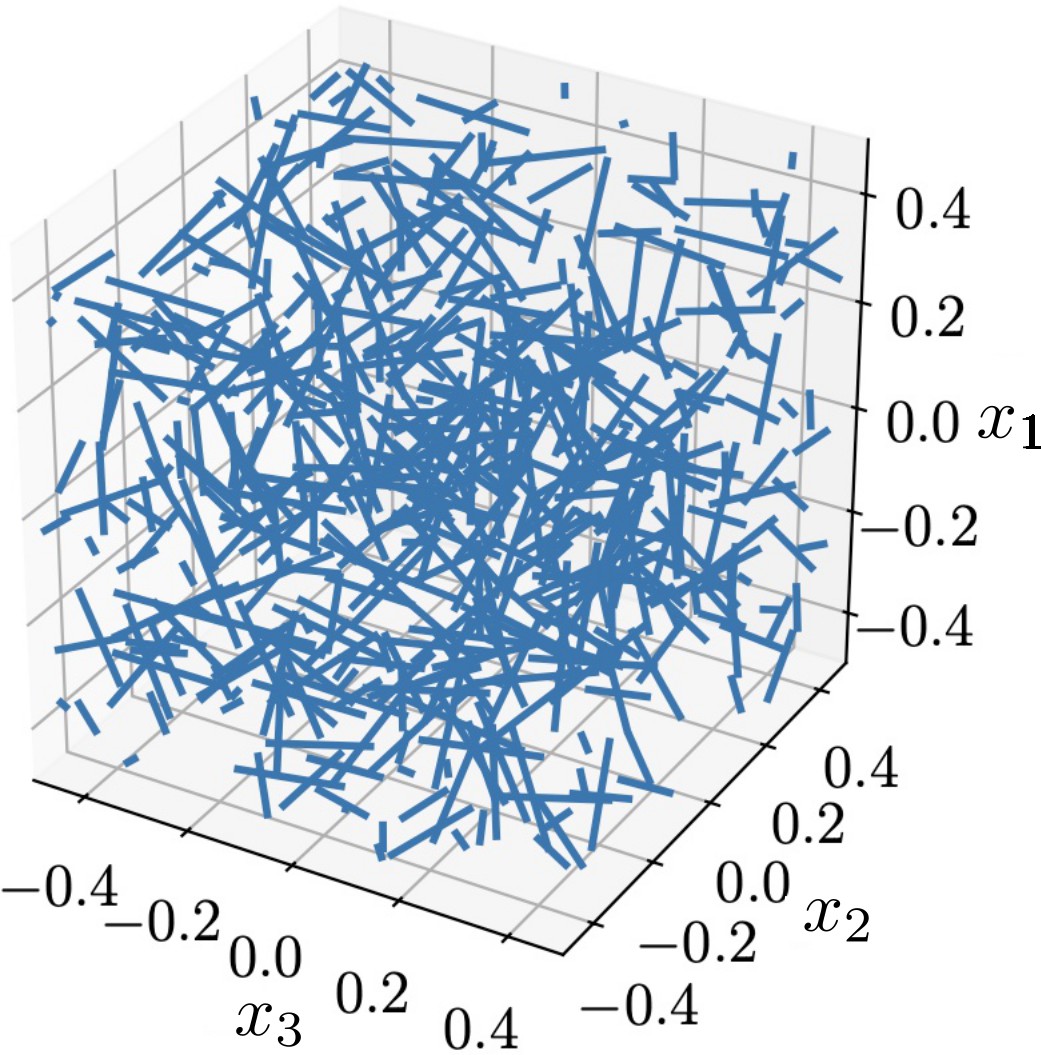}
    \caption{Unit cell with aligned fibers (left) and randomly oriented fibers (right), both with
      volume fraction $f_v=0.16$. Fibers have length $0.2$ and a circular cross section of radius $0.025$. The boundary conditions are periodic on all faces of the cube, resulting in fibers crossing between opposite faces.}
    \label{fig-inclusions-schemes}
\end{figure}

\begin{figure}[ht]
    \centering
       \includegraphics[width=0.45\linewidth]{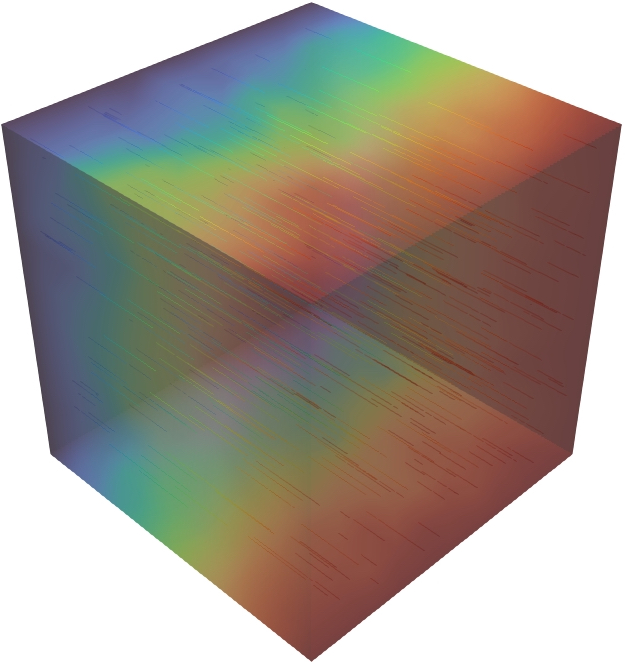}\hfill
        \includegraphics[width=0.45\linewidth]{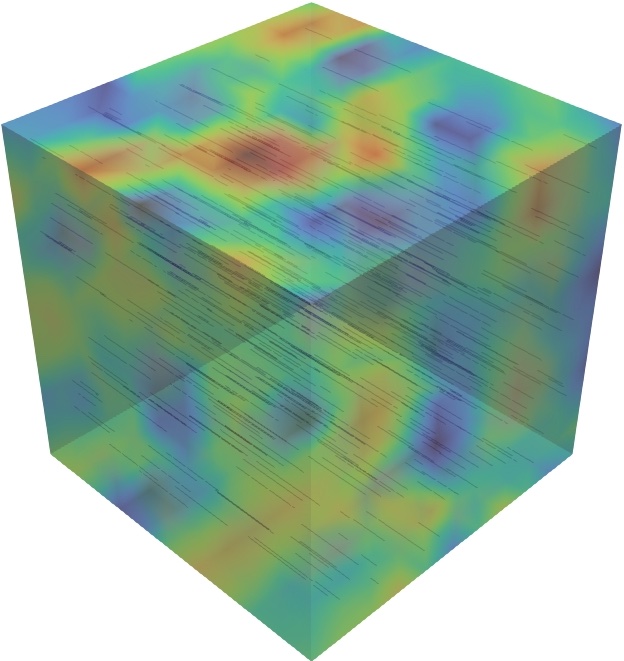}
   \vspace{5mm}
        \includegraphics[width=0.45\linewidth]{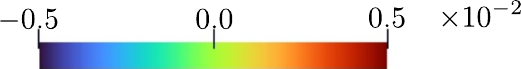}\hfill
        \includegraphics[width=0.45\linewidth]{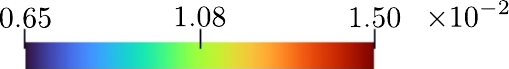}
    \caption{Displacement (left) and normal stress (right), both in the traction direction in the
      solid for the case of fibers with volume fraction $f_v=0.16$.}
    \label{fig-inclusions-oriented}
\end{figure}

\begin{figure}[ht]
    \centering
        \includegraphics[width=0.45\linewidth]{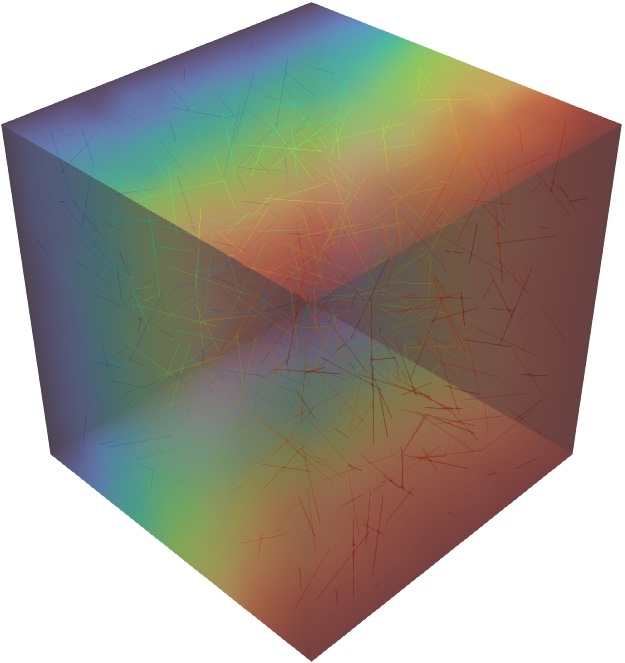}\hfill
        \includegraphics[width=0.45\linewidth]{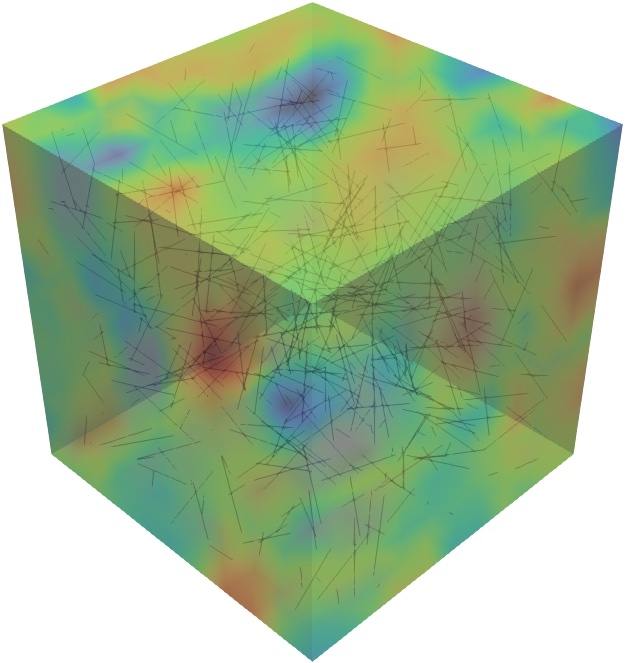}
   \vspace{5mm}
        \includegraphics[width=0.45\linewidth]{colorbar_inclusions_displacement.pdf}\hfill
        \includegraphics[width=0.45\linewidth]{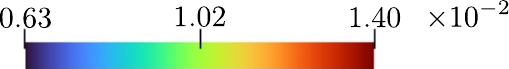}
    \caption{Displacement (left) and normal stress (right) in the traction direction in the solid.
      Fiber volume fraction $f_v=0.16$.}
    \label{fig-inclusions-random}
\end{figure}

The fibers have been randomly placed, and two configurations have been studied: fibers aligned with
the traction direction and fibers randomly oriented. Figure \ref{fig-inclusions-schemes} shows the
configurations for the two cases. The displacement
results for these cases can be seen in Figures~\ref{fig-inclusions-oriented} (aligned case)
and~\ref{fig-inclusions-random} (random orientation case). In both cases, displacement fluctuations
can be observed in Figures~\ref{fig-inclusions-oriented} and~\ref{fig-inclusions-random} as a consequence of the higher stiffness of the fibers. While in the case without fibers the normal
stress is homogeneous with a value of $\sigma_{33} = E_{\Omega} \; \bar{\varepsilon}_{33} = 0.01$, with
fiber reinforcement the stress component $\sigma_{33}$ takes values in the ranges $[0.0065, 0.0150]$ and $[0.0063, 0.0140]$ for the aligned and random cases, respectively.

We have studied the dependence of the effective Young's modulus of the cell on the fiber volume
fraction. For this purpose, six different volume fractions ($f_v = 0.005, 0.01, 0.02, 0.04, 0.08,
0.16$) were analyzed. To ensure statistical validity for each volume fraction, $20$ random
configurations were evaluated and then the mean and standard deviation of these tests were obtained.
The effective Young’s moduli are obtained obtained by calculating the average stress in the loading
direction in the solid and the fibers, $\bar{\sigma}_X$, as
\begin{equation}
\bar{\sigma}_{33} = \frac{1}{|\Omega|} \int_{\Omega} \sigma_{33}^{\Omega} \; \textrm{d}V + \sum_{i=1}^{N_{fibers}} \frac{1}{|\mathcal{S}_i|} \int_{\mathcal{S}_{i}} \sigma_{33}^{\mathcal{S}_i} \; \textrm{d} V\ ,
\end{equation}
where $\sigma_{33}^{\Omega}$ and $\sigma_{33}^{\mathcal{S}_i}$ are the normal stress in the traction direction in the solid $\Omega$ and structure $\mathcal{S}_i$, respectively. The effective Young’s modulus is, therefore,
\begin{equation}
\bar{E} = \frac{\bar{\sigma}_{33}}{\bar{\varepsilon}_{33}}\ .
\end{equation}

\begin{figure}[ht]
    \centering
    \includegraphics[width=0.9\textwidth]{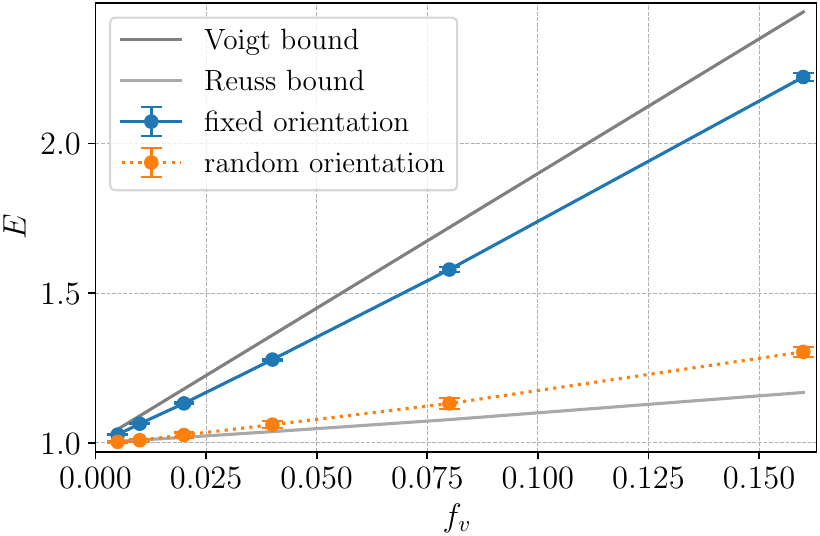}
    \caption{Effective Young’s modulus as a function of the fiber volume fraction. The points for
      the aligned case (blue) and the randomly oriented case (orange) are shown, along with the
      Voigt and Reuss bounds. Each point represents the average of 20 simulations with randomly
      located fibers, as well as their standard deviation.}
    \label{fig-inclusions-young}
\end{figure}

Figure \ref{fig-inclusions-young} shows the macroscopic Young’s modulus as a function of the fiber
volume fraction. Also the Voigt and Reuss bounds, $E_V = f_v E_{\mathcal{S}} + (1-f_v) E_{\Omega}$
and $E_R = \frac{1}{\frac{fv}{E_{\mathcal{S}}} + \frac{1-f_v}{E_{\Omega}}}$, are plotted. It can be
observed that in both cases two straight lines are defined within the bounds. In the case where the
fibers are aligned, the Young’s modulus obtained is close to the Voigt bound (around $90\%$). On the
other hand, the case with randomly oriented fibers lies closer to the Reuss bound, precisely as a
consequence of the fibers providing less stiffness due to their misalignment with respect to the traction direction.

\clearpage

\section{Summary of main results}\label{sec-summary}

In this work, we have presented a general framework for embedding structures into continuous solids
under small deformation conditions. By using an extension of the Arlequin method for coupling
solids, we have formulated a well-posed problem that also enables the straightforward definition of
stable mixed finite elements, as has been demonstrated. Remarkably, the method is completely
independent of both the solid and reduced model meshes, which can be arbitrary and incompatible.

An stability issue has been identified in the discrete case arising from the lack of coercivity under certain conditions. The problem depends on the structural model, the mesh sizes of the solid and the structure, their stiffnesses, and the interpolation order of the finite elements. The solution is simply to refine the structural mesh so that it is at least as fine as that of the solid, which has little impact on computational cost given that it is a reduced model. Moreover, this refinement ensures full integration in the new integrals introduced by the numerical method.

The numerical examples illustrate the properties of the proposed method. Given its generality, the method has proven effective in simulating all the motions of beams and shells, including cases where the base remains fixed but the fibers rotate, such as torsion in beams and pure shear in shells. All of this is achieved while significantly reducing the number of degrees of freedom, even enabling problems to be solved that would be unfeasible using a pure continuum-solid simulation, such as the case of micro inclusions.

Finally, let us note that the method proposed in this work is fairly general and can be readily
adapted to various structural models not addressed here (membranes, bars, plates, \dots), as well as to diffusion problems.
In particular, the extension to the finite strain follows the formalism outlined in the present
work, although some additional complications arise due to the potentially large size of the
displacements and rotations. These results will be published in a sequel.


\section*{Acknowledgements}
The authors would like to acknowledge Ansys, Inc. for the financial support and the enriching
discussions that significantly contributed to this work. This work has also received funding
from the Spanish Ministry of Science and Innovation under grant PID2021-128812OB-I00.

\clearpage

\appendix
\section{Coercivity analysis}
\label{app-instabilities}

We collect in this appendix some of the more technical bounds required to prove the coercivity of
the tied bodies, in both the continuum and the discrete cases.

\subsection{Continuum case}

For the embedding problem we can estimate the coercivity constant by identifying $\ker b$,
\begin{equation}
\label{eq:ker_b}
\ker{b} = \left\{ (u,\Sigma,\Theta) \in X \times Y \textrm{ s.t. }
  b(u,\Sigma,\Theta; \gamma,\mu) = 0,\
  \forall \; \combo{\gamma}{\mu} \in Y\right\}
\end{equation}
It is easy to see that $\ker{b}$ may be identified with $u |_{\tilde{S}}= \combo{\Sigma}{\Theta}$. In that case we have the following result,
\begin{equation}
\label{eq:coercivity_cont}
\begin{array}{lcl}
a(u,\Sigma,\Theta;u,\Sigma,\Theta)  & = & a_{LE}(u,u) + a_{ST}((\Sigma,\Theta),(\Sigma,\Theta)) \\
                                        & = & \int_{\Omega} \nabla u \; C_{\Omega}\; \nabla u  \dV + \int_{\mathcal{C}} \nabla \Sigma \; C_{\Sigma}\;  \nabla \Sigma + \nabla \theta \; C_{\theta} \; \nabla \theta \textrm{d}\mathcal{C} \\
                                        & \ge & \underline{C_{\Omega}} ||\nabla u||_{L^2(\Omega)}^2
                                       + \underline{C_{\Sigma}} ||\nabla \Sigma||_{L^2(\mathcal{C})}^2
                                           + \underline{C_{\theta}} ||\nabla \theta||_{L^2(\mathcal{C})}^2 \\
                                        & = & \epsilon_1 \underline{C_{\Omega}} ||\nabla u||_{L^2(\Omega)}^2
                                           + (1-\epsilon_1) \underline{C_{\Omega}} ||\nabla u||_{L^2(\Omega)}^2
                                           \\
                                        &   & + \underline{C_{\Sigma}} ||\nabla \Sigma||_{L^2(\mathcal{C})}^2
                                           + \underline{C_{\theta}} ||\nabla \theta||_{L^2(\mathcal{C})}^2 \\
                                        & \ge & \epsilon_1 \underline{C_{\Omega}} ||\nabla u||_{L^2(\Omega)}^2
                                           + \frac{(1-\epsilon_1) \underline{C_{\Omega}}}{C_{P_{\Omega}}} ||u||_{L^2(\Omega)}^2
                                           \\
                                        &   &
                                           + \underline{C_{\Sigma}} ||\nabla \Sigma||_{L^2(\mathcal{C})}^2
                                           + \underline{C_{\theta}} ||\nabla \theta||_{L^2(\mathcal{C})}^2  \\
                                        & = &
\epsilon_1 \underline{C_{\Omega}} ||\nabla u||_{L^2(\Omega)}^2
										  + \epsilon_2 \frac{(1-\epsilon_1) \underline{C_{\Omega}}}{C_{P_{\Omega}}} ||u||_{L^2(\Omega)}^2
										  \\
										&   &
										  + (1-\epsilon_2) \frac{(1-\epsilon_1) \underline{C_{\Omega}}}{C_{P_{\Omega}}} ||u||_{L^2(\Omega)}^2
										  \\
										&  &
                                           + \underline{C_{\Sigma}} ||\nabla \Sigma||_{L^2(\mathcal{C})}^2
                                           + \underline{C_{\theta}} ||\nabla \theta||_{L^2(\mathcal{C})}^2 \\
                                        & \ge &
\epsilon_1 \underline{C_{\Omega}} ||\nabla u||_{L^2(\Omega)}^2
										  + \epsilon_2 \frac{(1-\epsilon_1) \underline{C_{\Omega}}}{C_{P_{\Omega}}} ||u||_{L^2(\Omega)}^2
										 \\
									    &   &
										  + (1-\epsilon_2) \frac{(1-\epsilon_1) \underline{C_{\Omega}}}{C_{P_{\Omega}}} ||u||_{L^2(\tilde{S})}^2
										  \\
										&   &
                                           + \underline{C_{\Sigma}} ||\nabla \Sigma||_{L^2(\mathcal{C})}^2
                                           + \underline{C_{\theta}} ||\nabla \theta||_{L^2(\mathcal{C})}^2 \\
                                        & = &
\epsilon_1 \underline{C_{\Omega}} ||\nabla u||_{L^2(\Omega)}^2
										  + \epsilon_2 \frac{(1-\epsilon_1) \underline{C_{\Omega}}}{C_{P_{\Omega}}} ||u||_{L^2(\Omega)}^2
										  \\
									   &   &
										  + (1-\epsilon_2) \frac{(1-\epsilon_1) \underline{C_{\Omega}}}{C_{P_{\Omega}}} ||\combo{\Sigma}{\Theta}||_{L^2(\tilde{S})}^2
										  \\
									   &   &
                                           + \underline{C_{\Sigma}} ||\nabla \Sigma||_{L^2(\mathcal{C})}^2
                                           + \underline{C_{\theta}} ||\nabla \theta||_{L^2(\mathcal{C})}^2 \\
                                        & \ge &
\epsilon_1 \underline{C_{\Omega}} ||\nabla u||
   									  + \epsilon_2 \frac{(1-\epsilon_1) \underline{C_{\Omega}}}{C_{P_{\Omega}}} ||u||_{L^2(\Omega)}^2
   									      \\
   									    &   &
										  + (1-\epsilon_2) \frac{(1-\epsilon_1) \underline{C_{\Omega}}}{C_{P_{\Omega}}} \left( \alpha_{\Sigma} ||\Sigma||_{L^2(\mathcal{C})}^2 + \alpha_{\theta}  ||\theta||_{L^2(\mathcal{C})}^2 \right)
										  \\
										&   &
                                           + \underline{C_{\Sigma}} ||\nabla \Sigma||_{L^2(\mathcal{C})}^2
                                           + \underline{C_{\theta}} ||\nabla \theta||_{L^2(\mathcal{C})}^2 \\
                                        & \ge &
 \min \left(\epsilon_1, \epsilon_2 \frac{(1-\epsilon_1)}{C_{P_{\Omega}}}\right) \underline{C_{\Omega}} ||u||_{H^1(\Omega)}
 									      \\
 									    &   &
                                           + \min \left( (1-\epsilon_2)\frac{(1-\epsilon_1)\underline{C_{\Omega}}}{C_{P_{\Omega}}} \alpha_{\Sigma}, \underline{C_{\Sigma}}\right) ||\Sigma||_{H^1(\mathcal{C})}^2
                                           \\
                                         &   &
                                           + \min \left( (1-\epsilon_2)\frac{(1-\epsilon_1)\underline{C_{\Omega}}}{C_{P_{\Omega}}} \alpha_{\theta}, \underline{C_{\theta}}\right) ||\theta||_{H^1(\mathcal{C})}^2\ ,
\end{array}
\end{equation}
where $\underline{C_i}$ are lower bounds of the corresponding elasticity tensor, $C_{P_{\Omega}}$ is
Poincare's constant and $\alpha_{\Sigma}$, $\alpha_{\theta} > 0$ are constants that depend on the
fiber's geometry. Taking $0<\epsilon_1,\epsilon_2<1$ we obtain lower bounds for the coercivity constant.

\subsection{Coercivity analysis for structures with only translations}
For the simple case in which the structure only has translational degrees of freedom, the discrete structure bilinear form reduces to
\begin{equation}
a_{ST}(\Sigma_{\hb},\Gamma_{\hb}) = \int_{\mathcal{C}} D \Sigma_{\hb} \cdot C_{\Sigma} \; D\Gamma_{\hb} \textrm{d}\mathcal{C}\ ,
\end{equation}
and the coupling term between solid and structure becomes
\begin{equation}
b(\lambda_{\hb},\Gamma_{\hb}-v_h)= \int_{S} \left[ \lambda_{\hb} \cdot (\Gamma_{\hb}-v_h) + \ell^2 D\lambda_{\hb} \cdot D(\Gamma_{\hb} - v_{h} )\right] \dV\ .
\end{equation}

We distinguish here two cases depending on the mesh size of the structure respect to the solid.

\subsubsection{Finer structure mesh}
For only translation degrees of freedom and finer structure mesh, we have $\tilde{\Omega}_h
\subseteq \tilde{C}_h$, where $\tilde{\Omega}_h$ is the part of the discrete solid which coincides
with the discrete structure, and then $u_h|_{\tilde{\Omega}_h} \in \tilde{C}_h$. Taking this into
account,  analogous expressions as in Eq.~\eqref{eq:coercivity_cont} can be obtained by neglecting
the rotation terms for the discrete case.
As in the continuum case, the structure's lack of coercivity is compensated by the solid's.

\subsubsection{Coarser structure mesh}
The case of the coarser mesh of the structure is substantially different from the previous one. In this case, $\ker b$ is no longer identified with solutions that are equal in the common region for the solid displacements and the lifting of the structure. Instead, the kernel is identified with the lifting of the structure being equal to the projection of the solid displacements onto a space with the mesh size of the structure. Designating $\uhb$ as the projection of the solid displacement $u_h$ onto the functional space given by $\hb$, it can be seen that in $\ker b_{\hb}$, $\Gamma_{\hb} = \vhb$. Taking $\lambda_{\hb} = \Gamma_{\hb} - \vhb$, considering that the projection $\vhb$ is orthogonal in the finite element space, and that the space of the structure is contained within that of the solid,

\begin{equation}
\label{eq:kernel_coarser}
\begin{array}{lcl}
0 & =&  \int_{\tilde{S}_h} \left[ (\Gamma_{\hb} - \vhb) \cdot (\Gamma_{\hb}-v_{\hs}) + \ell^2 D(\Gamma_{\hb} - \vhb) \cdot D(\Gamma_{\hb} - v_{\hs} )\right] \dV \\
 & = & \int_{\tilde{S}_h} \left[ (\Gamma_{\hb} - \vhb) \cdot (\Gamma_{\hb}-v_{\hs} + \vhb -\vhb )\right] \dV \\
 & & +  \int_{\tilde{S}_h} \ell^2 \left[ D(\Gamma_{\hb} - \vhb) \cdot D(\Gamma_{\hb} - v_{\hs} +\vhb - \vhb )\right] \dV \\
 & = & \int_{\tilde{S}_h} \left[ (\Gamma_{\hb} - \vhb) \cdot (\Gamma_{\hb}-\vhb ) + \ell^2 D(\Gamma_{\hb} - \vhb) \cdot D(\Gamma_{\hb} - \vhb )\right] \dV \\
  &  & + \int_{\tilde{S}_h} \left[ (\Gamma_{\hb} - \vhb) \cdot (\vhb - v_{\hs} ) + \ell^2 D(\Gamma_{\hb} - \vhb) \cdot D(\vhb - v_{\hs} )\right] \dV \\
  & = & || \Gamma_{\hb} - \vhb ||_{H^1(\tilde{S}_h)}^2\ ,
\end{array}
\end{equation}
and we get $\Gamma_{\hb} = \vhb$ in $\ker b_h$.

Taking into account~\eqref{eq:kernel_coarser}, Poincare's inequality and the approximation estimates, we can approximate the coercivity constant,
\begin{equation}
\label{eq:coercivity_disc}
\begin{array}{lcl}
a(u_{\hs},\Sigma_{\hb}; u_{\hs},\Sigma_{\hb})  & = & a_{LE}(u_{\hs},u_{\hs}) + a_{ST}((\Sigma_{\hb},\Sigma_{\hb}) \\
                                        & = & \int_{\Omega_{\hs}} D u_{\hs} \; C_{\Omega}\; D u_{\hs}  \dV + \int_{\mathcal{C}} D \Sigma_{\hb} \; C_{\Sigma}\;  D \Sigma_{\hb} \textrm{d}\mathcal{C} \\

                                        & \ge & \underline{C_{\Omega}} ||D u_{\hs}||_{L^2(\Omega_{\hs})}^2
                                       + \underline{C_{\Sigma}} ||D \Sigma_{\hb}||_{L^2(\mathcal{C}_{\hb})}^2
                                       \\
                                       & = & \epsilon_1 \underline{C_{\Omega}} ||D u_{\hs}||_{L^2(\Omega_h)}^2
                                           + (1-\epsilon_1) \underline{C_{\Omega}} ||D (u_{\hs})||_{L^2(\Omega_{\hs})}^2 \\
                                           & & + \underline{C_{\Sigma}} ||D \Sigma_{\hb}||_{L^2(\mathcal{C}_{\hb})}^2
                                           \\
                                       & \ge & \epsilon_1 \underline{C_{\Omega}} ||D u_{\hs}||_{L^2(\Omega_h)}^2
                                           + \frac{(1-\epsilon_1) \underline{C_{\Omega}}}{C_{P_{\Omega_{\hs}}}} ||u_{\hs}||_{L^2(\Omega_{\hs})}^2 \\
                                           & & + \underline{C_{\Sigma}} ||D \Sigma_{\hb}||_{L^2(\mathcal{C}_{\hb})}^2
                                           \\
                                            & = & \epsilon_1 \underline{C_{\Omega}} ||D u_{\hs}||_{L^2(\Omega_h)}^2
                                           + \epsilon_2  \frac{(1-\epsilon_1) \underline{C_{\Omega}}}{C_{P_{\Omega_{\hs}}}} ||u_{\hs}||_{L^2(\Omega_{\hs})}^2 \\
                                           & & + (1-\epsilon_2)  \frac{(1-\epsilon_1) \underline{C_{\Omega}}}{C_{P_{\Omega_{\hs}}}} ||u_{\hs}||_{L^2(\Omega_{\hs})}^2
                                           + \underline{C_{\Sigma}} ||D \Sigma_{\hb}||_{L^2(\mathcal{C}_{\hb})}^2

                                           \\
                                           & \ge &  \epsilon_1 \underline{C_{\Omega}} ||D u_{\hs}||_{L^2(\Omega_h)}^2
                                           + \epsilon_2  \frac{(1-\epsilon_1) \underline{C_{\Omega}}}{C_{P_{\Omega_{\hs}}}} ||u_{\hs}||_{L^2(\Omega_{\hs})}^2 \\
                                           & & + (1-\epsilon_2)  \frac{(1-\epsilon_1) \underline{C_{\Omega}}}{C_{P_{\Omega_{\hs}}}} ||u_{\hs}||_{L^2(\tilde{S}_{h})}^2
                                            + \underline{C_{\Sigma}} ||D \Sigma_{\hb}||_{L^2(\mathcal{C}_{\hb})}^2
                                           \\
                                            & = &  \epsilon_1 \underline{C_{\Omega}} ||D u_{\hs}||_{L^2(\Omega_h)}^2
                                           + \epsilon_2  \frac{(1-\epsilon_1) \underline{C_{\Omega}}}{C_{P_{\Omega_{\hs}}}} ||u_{\hs}||_{L^2(\Omega_{\hs})}^2 \\
                                           & & + (1-\epsilon_2)  \frac{(1-\epsilon_1) \underline{C_{\Omega}}}{C_{P_{\Omega_{\hs}}}}||u_{\hs}+\uhb-\uhb||_{L^2(\tilde{S}_{h})}^2 \\
                                           & & + \underline{C_{\Sigma}} ||D \Sigma_{\hb}||_{L^2(\mathcal{C}_{\hb})}^2
                                           \\
                                            & \ge &  \epsilon_1 \underline{C_{\Omega}} ||D u_{\hs}||_{L^2(\Omega_h)}^2
                                           + \epsilon_2  \frac{(1-\epsilon_1) \underline{C_{\Omega}}}{C_{P_{\Omega_{\hs}}}} ||u_{\hs}||_{L^2(\Omega_{\hs})}^2 \\
                                           & & + (1-\epsilon_2)  \frac{(1-\epsilon_1) \underline{C_{\Omega}}}{C_{P_{\Omega_{\hs}}}} \left(||\uhb||_{L^2(\tilde{S}_{h})}^2 -||u_{\hs} - \uhb||_{L^2(\tilde{S}_{h})}^2\right) \\
                                         & &  + \underline{C_{\Sigma}} ||D \Sigma_{\hb}||_{L^2(\mathcal{C}_{\hb})}^2
                                         \\
                                            & \ge &  \epsilon_1 \underline{C_{\Omega}} ||D u_{\hs}||_{L^2(\Omega_h)}^2
                                           + \epsilon_2  \frac{(1-\epsilon_1) \underline{C_{\Omega}}}{C_{P_{\Omega_{\hs}}}} ||u_{\hs}||_{L^2(\Omega_{\hs})}^2 \\
                                           & &  -(1-\epsilon_2)  \frac{(1-\epsilon_1) \underline{C_{\Omega}}}{C_{P_{\Omega_{\hs}}}} \beta \hb^{2(p+1)}|u_{\hs}|_{H^{p+1}(\tilde{S}_{h})}^2 \\
                                         & &  + (1-\epsilon_2)  \frac{(1-\epsilon_1) \underline{C_{\Omega}}}{C_{P_{\Omega_{\hs}}}} ||\uhb||_{L^2(\tilde{S}_h)}^2
                                         + \underline{C_{\Sigma}} ||D \Sigma_{\hb}||_{L^2(\mathcal{C}_{\hb})}^2
                                         \\
                                            & = &  \epsilon_1 \underline{C_{\Omega}} ||D u_{\hs}||_{L^2(\Omega_h)}^2
                                           + \epsilon_2  \frac{(1-\epsilon_1) \underline{C_{\Omega}}}{C_{P_{\Omega_{\hs}}}} ||u_{\hs}||_{L^2(\Omega_{\hs})}^2 \\
                                           & &  -(1-\epsilon_2)  \frac{(1-\epsilon_1) \underline{C_{\Omega}}}{C_{P_{\Omega_{\hs}}}} \beta \hb^{2(p+1)}|u_{\hs}|_{H^{p+1}(\tilde{S}_{h})}^2 \\
                                         & &  +(1-\epsilon_2)  \frac{(1-\epsilon_1) \underline{C_{\Omega}}}{C_{P_{\Omega_{\hs}}}} ||\Sigma_{\hb}||_{L^2(\tilde{S}_{h})}^2
                                         + \underline{C_{\Sigma}} ||D \Sigma_{\hb}||_{L^2(\mathcal{C}_{\hb})}^2
                                         \\
                                            & \ge & \underline{C_{\Sigma}}
                                            \left(
                                            \min \left(\frac{\epsilon_1 \underline{C_{\Omega}}}{\underline{C_{\Sigma}}}, \epsilon_2  \frac{(1-\epsilon_1) \underline{C_{\Omega}}}{C_{P_{\Omega_{\hs} \; \underline{C_{\Sigma}}}}} \right) ||u_{\hs}||_{H^1(\Omega_h)}^2  \right.
                                            \\
                                           & &  -(1-\epsilon_2)  \frac{(1-\epsilon_1) \underline{C_{\Omega}}}{C_{P_{\Omega_{\hs}}} \; \underline{C_{\Sigma}}} \beta \hb^{2(p+1)}|u_{\hs}|_{H^{p+1}(\Omega_{\hs})}^2 \\
                                         & & \left.
                                         + \min \left( (1-\epsilon_2)  \frac{(1-\epsilon_1) \; |\fiber| \; \underline{C_{\Omega}}}{C_{P_{\Omega_{\hs}}} \; \underline{C_{\Sigma}}}, 1 \right)  ||\Sigma_{\hb}||_{H^1(\mathcal{C}_{\hb})}^2
                                          \right)\ ,
\end{array}
\end{equation}
for $0 < \epsilon_1, \epsilon_2 < 1$ and $\beta$ is a constant that does not depend on the mesh. Finally, $p$ is the order of the finite element approximation.


\bibliographystyle{\mybibstyle}
\bibliography{driver.bib}

\end{document}